\documentclass[aos,MSNbibl,seceqn,citesort,dvips]{arximspdf}
\usepackage{graphicx}
\doi{10.1214/11-AOS920}
\volume{39}
\issue{5}
\pubyear{2011}
\firstpage{2626}
\lastpage{2657}

\makeatletter

\newtheorem{theorem}{Theorem}[section]
\newtheorem{lemma}{Lemma}[section]
\newtheorem{proposition}{Proposition}[section]

\newproclaim{assumption}[theorem]{Assumption}

\renewcommand{\a}{{\alpha}}
\renewcommand{\b}{{\beta}}
\renewcommand{\d}{{\delta}}
\newcommand{\g}{{\gamma}}
\newcommand{\h}{{\eta}}
\renewcommand{\k}{{\kappa}}
\renewcommand{\l}{{\lambda}}
\newcommand{\m}{{\mu}}
\newcommand{\n}{{\nu}}
\newcommand{\NN}{\mathbb{N}}
\renewcommand{\r}{{\rho}}
\newcommand{\RR}{\mathbb{R}}
\renewcommand{\t}{{\tau}}

\newcommand{\ra}{\rightarrow}
\newcommand{\da}{\downarrow}
\newcommand{\ua}{\uparrow}

\newcommand{\E}{\mathrm{E}}
\renewcommand{\Pr}{\mathrm{P}}
\newcommand{\var}{\operatorname{var}}
\newcommand{\sdev}{\operatorname{sd}}
\newcommand{\Cov}{\operatorname{Cov}}
\newcommand{\given}{|}

\newcommand{\weak}{\rightsquigarrow}
\newcommand{\tr}{\operatorname{tr}}
\newcommand{\RV}{\mathcal{R}}
\newcommand{\eps}{\varepsilon}
\newcommand{\wt}{\tilde}
\renewcommand{\S}{\mathcal{S}}

\makeatother

\begin{document}
\begin{frontmatter}

\title{Bayesian inverse problems with Gaussian priors}
\runtitle{Bayesian inverse problems with Gaussian priors}

\begin{aug}
\author[A]{\fnms{B. T.} \snm{Knapik}\corref{}\thanksref{t1}\ead[label=e1]{b.t.knapik@vu.nl}},
\author[A]{\fnms{A. W.} \snm{van der Vaart}\ead[label=e2]{aad@cs.vu.nl}}
\and
\author[B]{\fnms{J. H.} \snm{van Zanten}\ead[label=e3]{j.h.v.zanten@tue.nl}}
\runauthor{B. T. Knapik, A. W. van der Vaart and J. H. van Zanten}
\affiliation{VU University Amsterdam, VU University Amsterdam and\break
Eindhoven University of Technology}
\address[A]{B. T. Knapik\\
A. W. van der Vaart\\
Department of Mathematics\\
Faculty of Sciences\\
VU University Amsterdam\\
De Boelelaan 1081a \\
1081 HV Amsterdam\\
The Netherlands\\
\printead{e1}\\
\hphantom{E-mail: }\printead*{e2}}
\address[B]{J. H. van Zanten\\
Department of Mathematics\\
Eindhoven University of Technology\\
Eindhoven\\
The Netherlands\\
\printead{e3}} %adresu isvedimo komanda gale!
\end{aug}

\thankstext{t1}{Supported in part by the Netherlands Organization for
Scientific Research NWO.}

% HISTORY:
\received{\smonth{3} \syear{2011}}
\revised{\smonth{8} \syear{2011}}

% ABSTRACT
%
\begin{abstract}
The posterior distribution in a nonparametric inverse problem is shown
to contract to the true parameter at a rate that depends on the
smoothness of the parameter, and the smoothness and scale of the prior.
Correct combinations of these characteristics lead to the minimax rate.
The frequentist coverage of credible sets is shown to depend on the
combination of prior and true parameter, with smoother priors leading
to zero coverage and rougher priors to conservative coverage. In the
latter case credible sets are of the correct order of magnitude. The
results are numerically illustrated by the problem of recovering a
function from observation of a noisy version of its primitive.
\end{abstract}

% KEYWORDS
%

\end{frontmatter}

%s1 ###
%se1 #&#
\section{Introduction}\label{sec1}

In this paper we study a Bayesian approach to estimating
a parameter $\m$ from an observation $Y$ following the model
%
%e1.1 ###
%e1.1 #&#
%
\begin{equation}
\label{EqProblem}
Y = K\m+ \frac{1}{\sqrt{n}}Z.
\end{equation}
The unknown parameter $\m$ is an element
of a separable Hilbert space $H_1$, and is mapped into another
Hilbert space $H_2$ by a known, injective, continuous
linear operator $K\dvtx H_1 \to H_2$.
The image $K\m$ is perturbed by unobserved, scaled Gaussian white
noise $Z$.
There are many special examples of this infinite-dimensional regression model,
which can also be viewed as an idealized
version of other statistical models, including density estimation.
The inverse problem of estimating $\m$ has
been studied by both statisticians and numerical
mathematicians (see, e.g., \cite
{Donoho,Cavalier,Munk,RuymgaartII,RuymgaartIII,Stuart}
for reviews), but rarely from a theoretical Bayesian perspective;
exceptions are \cite{Cox} and \cite{Simoni}.

The Bayesian approach to (\ref{EqProblem}) consists of putting a prior
on the parameter~$\m$, and computing the posterior distribution. We
study Gaussian priors, which are conjugate to the model, so that the
posterior distribution is also Gaussian and easy to derive. Our interest
is in studying the properties of this posterior distribution, under
the ``frequentist'' assumption that the data $Y$ has been generated
according to the model (\ref{EqProblem}) with a~given ``true''
parameter $\m_0$. We investigate whether and at what rate the
posterior distributions contract to $\m_0$ as $n\ra\infty$ (as in
\cite{GGvdV}), but have as main interest the performance of
credible sets for measuring the uncertainty about the parameter.

A Bayesian \textit{credible set} is defined as a central region in the
posterior distribution of specified posterior probability, for
instance, 95\%. As a consequence of the Bernstein--von Mises theorem
credible sets for smooth \textit{finite-dimensional}
parametric models are asymptotically equivalent to confidence regions
based on the maximum likelihood estimator (see, e.g., \cite{vdVAS},
Chapter~10), under mild conditions on the prior. Thus, ``Bayesian
uncertainty'' is equivalent to ``frequentist uncertainty'' in these
cases, at least for large~$n$. However, there is no corresponding
Bernstein--von Mises theorem in nonparametric
Bayesian inference, as noted in \cite{Freedman}. The
performance of Bayesian credible sets in these situations has
received little attention, although in practice such sets are
typically provided as indicators of uncertainty, for instance, based
on the spread of the output of a (converged) MCMC run. The paper \cite{Cox}
did tackle this issue and came to the alarming conclusion
that Bayesian credible sets have frequentist coverage zero.
If this were true, many data analysts would
(justifiably) distrust the spread in the posterior distribution as a
measure of uncertainty. For other results see
\cite{Bontemps,Ghosal99,Ghosal00} and~\cite{Leahu}.

The model considered in \cite{Cox} is equivalent to
our model (\ref{EqProblem}), and a good starting point
for studying these issues. More precisely, the conclusion
of~\cite{Cox} is that \textit{for almost every parameter $\m_0$ from
the prior the
coverage of a credible set} (\textit{of any level}) \textit{is} 0. In the present paper
we show that this is only part of the story, and, taken by itself,
the conclusion is misleading. The coverage
depends on the true parameter $\m_0$ and the prior together,
and it can be understood in terms of a bias-variance trade-off,
much as the coverage of frequentist nonparametric procedures.
A nonparametric procedure that oversmoothes the truth (too big a bandwidth
in a~frequentist procedure, or a prior that puts too much weight
on ``smooth'' parameters) will be biased, and a confidence
or credible region based on such a procedure will be both
too concentrated and wrongly located, giving zero coverage. On the other
hand, undersmoothing does work (to a certain extent),
also in the Bayesian setup, as we show below.
In this light we reinterpret the conclusion
of \cite{Cox} to be valid only in the oversmoothed case (notwithstanding
a conjecture to the contrary in the Introduction of this paper;
see page 905, answer to objection~4). In the undersmoothed case
credible regions are conservative in general, with coverage
tending to 1. The good news is that typically they are
of the correct order of magnitude, so that they do give
a reasonable idea of the uncertainty in the estimate.

Of course, whether a prior under- or oversmoothes depends
on the regularity of~the true parameter. In practice, we may not
want to consider this known, and adapt the prior smoothness
to the data. In this paper we do consider the effect of changing
the ``length scale'' of a prior, but do not study data-dependent
length scales. The effect of setting the latter by, for example, an
empirical or full Bayes method will require further study.

Credible sets are by definition ``central regions'' in the
posterior distribution. Because the posterior distribution is a random
probability measure on the Hilbert space $H_1$, a ``central ball'' is
a natural shape of such a set, but it has the disadvantage that it is
difficult to visualize. If the Hilbert space is a function space,
then \textit{credible bands} are more natural. These correspond to
simultaneous credible intervals for the function at a point, and
can be obtained from the (marginal) posterior distributions of a set
of linear functionals. Besides the full posterior distribution, we
therefore study its marginals for linear functionals. The same issue
of the dependence of coverage on under- and oversmoothing arises, except
that ``very smooth'' linear functionals cancel the inverse nature
of the problem, and do allow a~Bernstein--von Mises theorem for a
large set of priors. Unfortunately point evaluations are usually
not smooth in this sense.

Thus, we study two aspects of inverse problems---recovering the full
parameter $\m$ (Section \ref{SectionFull}) and recovering linear
functionals (Section \ref{SectionLinear}). We obtain the rate of
contraction of the posterior distribution in both settings, in its
dependence on parameters of the prior. Furthermore, and most
importantly, we study the ``frequentist'' coverage of credible regions
for $\m$ in both settings, for the same set of priors. In the next
section we give a more precise statement of the problem, and in
Section~\ref{SectionPriorPosterior} we describe the priors that we
consider and derive the corresponding posterior distributions. In
Section \ref{SectionVolterra} we illustrate the results by
simulations and pictures in the particular example that $K$ is the
Volterra operator. Technical proofs are placed in
Sections~\ref{SectionProofs} and \ref{SectionTechnical}
at the end of the paper.

Throughout the paper $\langle\cdot,\cdot\rangle_1$ and \mbox{$\|\cdot\|_1$},
and $\langle\cdot,\cdot\rangle_2$ and \mbox{$\|\cdot\|_2$} denote the inner
products and norms of the Hilbert spaces $H_1$ and $H_2$. The adjoint
of an operator $A$ between two Hilbert spaces is denoted by $A^T$.
The Sobolev space $S^\b$ with its norm \mbox{$\|\cdot\|_\b$} is defined
in (\ref{EqDefSobolev}). For two sequences $(a_n)$ and~$(b_n)$ of numbers
$a_n\asymp b_n$ means that $|a_n/b_n|$ is bounded away from zero
and infinity as $n\ra\infty$, and $a_n\lesssim b_n$ means that
$a_n/b_n$ is bounded.

%s2 ###
%se2 #&#
\section{Detailed description of the problem}
\label{SectionDetailed}
The noise process $Z$ in (\ref{EqProblem}) is the
standard normal or \textit{iso-Gaussian process}
for the Hilbert space $H_2$. Because this is
not realizable as a random element in $H_2$, the model (\ref{EqProblem})
is interpreted in process form (as in \cite{Munk}). The iso-Gaussian process
is the zero-mean Gaussian process $Z=(Z_h\dvtx h\in H_2)$ with covariance function
$\E Z_{h}Z_{h'}=\langle h,h'\rangle_2$, and the measurement
equation (\ref{EqProblem}) is interpreted in that we
observe a Gaussian process $Y=(Y_h\dvtx h\in H_2)$ with
mean and covariance functions
%
%e2.1 ###
%e2.1 #&#
%
\begin{equation}
\label{EqMeanCovY}
\E Y_h=\langle K\m,h\rangle_2,\qquad
\operatorname{cov} (Y_h,Y_{h'})=\frac1n\langle h,h'\rangle_2.
\end{equation}
Sufficiency considerations show that it is statistically equivalent
to observe the subprocess $(Y_{h_i}\dvtx i\in\NN)$, for any orthonormal
basis $h_1,h_2,\ldots$ of $H_2$.

If the operator $K$ is compact, then the \textit{spectral decomposition}
of the self-adjoint operator $K^TK\dvtx H_1\to H_1$ provides a convenient
basis. In the compact case the operator $K^TK$ possesses countably
many positive eigenvalues $\k_i^2$ and there is a corresponding
orthonormal basis $(e_i)$ of $H_1$ of eigenfunctions (hence,
$K^TKe_i=\k_i^2e_i$ for $i\in\NN$; see, e.g., \cite{Rudin}). The
sequence $(f_i)$ defined by $Ke_i=\k_if_i$ forms an orthonormal
``conjugate'' basis of the range of $K$ in $H_2$. An element $\m\in
H_1$ can be identified with its sequence $(\m_i)$ of coordinates
relative to the eigenbasis $(e_i)$, and its image
$K\m=\sum_i\m_iKe_i=\sum_i\mu_i\k_if_i$ can be identified with its
coordinates $(\m_i\k_i)$ relative to the conjugate basis~$(f_i)$. If
we write $Y_i$ for $Y_{f_i}$, then (\ref{EqMeanCovY}) shows that
$Y_1,Y_2,\ldots$ are independent Gaussian variables with means $\E
Y_i=\mu_i\k_i$ and variance $1/n$. Therefore, a~concrete equivalent
description of the statistical problem is to \textit{recover the
sequence $(\m_i)$ from independent observations $Y_1,Y_2,\ldots$
with $N(\mu_i\k_i,1/n)$-distributions}.

In the following we do not require $K$ to be compact, but
we do assume the existence of an orthonormal basis of eigenfunctions of $K^TK$.
The main additional example we then cover is the \textit{white noise
model}, in which $K$ is the identity operator. The description of
the problem remains the same.

If $\k_i\ra0$, this problem is ill-posed, and the recovery of
$\m$ from $Y$ an \textit{inverse problem}. The ill-posedness can be
quantified by the speed of decay $\k_i\da0$. Although the tools are
more widely applicable, we limit ourselves to the \textit{mildly
ill-posed} problem (in the terminology of \cite{Cavalier}) and
assume that the decay is polynomial: for some $p\ge0$,
\[
\k_i \asymp i^{-p}.
\]
Estimation of $\m$ is harder if the decay is faster (i.e., $p$ is larger).

The difficulty of estimation may be measured by the minimax
risks over the scale of \textit{Sobolev spaces} relative to the
orthonormal basis $(e_i)$ of eigenfunctions of $K^TK$. For $\b>0$ define
%
%e2.2 ###
%e2.2 #&#
%
\begin{equation}
\label{EqDefSobolev}
\|\m\|_\b= \sqrt{\sum_{i=1}^\infty\m_i^2i^{2\b}}
\qquad\mbox{if }
\m= \sum_{i=1}^\infty\m_ie_i.
\end{equation}
Then the Sobolev space of order $\b$ is
$S^\b= \{\m\in H_1\dvtx \|\m\|_\b<\infty\}$.
The minimax rate of estimation over the unit\vadjust{\goodbreak} ball of this space
relative to the loss $\|t-\m\|_1$ of an estimate $t$ for $\m$ is
$n^{-\b/(1+2\b+2p)}$.
This rate is attained by various ``regularization'' methods,
such as generalized \textit{Tikhonov} and \textit{Moore--Penrose}
regularization \cite{Mair,Bertero,Cavalier,Goldenshluger,Munk}.
The Bayesian approach is closely connected to these methods:
in Section~\ref{SectionPriorPosterior}
the posterior mean is shown to be a regularized estimator.

Besides recovery of the full parameter $\m$, we consider estimating
linear functionals $L\m$. The minimax rate for such functionals over
Sobolev balls depends on $L$ as well as on the parameter of the Sobolev
space. Decay of the coefficients of $L$ in the eigenbasis may alleviate
the level of ill-posedness, with rapid decay even bringing the
functional in the domain of ``regular'' $n^{-1/2}$-rate estimation.

%s3 ###
%se3 #&#
\section{Prior and posterior distributions}
\label{SectionPriorPosterior}
We assume a mean-zero Gaussian prior for the parameter $\m$.
In the next three paragraphs we recall some essential facts on
Gaussian distributions on Hilbert spaces.

A \textit{Gaussian distribution} $N(\n,\Lambda)$
on the Borel sets of the Hilbert space~%
$H_1$ is characterized by a \textit{mean} $\n$, which
can be any element of $H_1$, and a \textit{covariance operator}
$\Lambda\dvtx H_1\to H_1$, which is a nonnegative-definite,
self-adjoint, linear operator
\textit{of trace class}: a compact operator with
eigenvalues~$(\l_i)$ that are summable $\sum_{i=1}^\infty\l
_i<\infty$
(see, e.g., \cite{Skorohod}, pages 18--20). A random element $G$ in
$H_1$ is
$N(\n,\Lambda)$-distributed if and only if
the stochastic process $(\langle G, h\rangle_1\dvtx h\in H_1)$
is a Gaussian process with mean and covariance functions
%
%e3.1 ###
%e3.1 #&#
%
\begin{equation}
\label{EqProcessG}
\E\langle G, h\rangle_1=\langle\n, h\rangle_1,\qquad
\operatorname{cov} (\langle G, h\rangle_1,\langle G,
h'\rangle_1)
=\langle h, \Lambda h'\rangle_1.
\end{equation}
The coefficients $G_i=\langle G,\varphi_i\rangle_1$
of $G$ relative to an orthonormal eigenbasis $(\varphi_i)$
of $\Lambda$ are independent, univariate Gaussians
with means the coordinates $(\n_i)$ of the mean vector $\n$
and variances the eigenvalues $\l_i$.

The iso-Gaussian process $Z$ in (\ref{EqProblem}) may be thought of as
a $N(0,I)$-distributed Gaussian element, for $I$ the identity operator
(on $H_2$), but as $I$ is not of trace class, this distribution is not
realizable as a proper random element in $H_2$. Similarly, the data
$Y$ in (\ref{EqProblem}) can be described as having a~$N(K\m,
n^{-1}I)$-distribution.

For a stochastic process $W=(W_h\dvtx h\in H_2)$ and a
continuous, linear operator $A\dvtx H_2\to H_1$, we define the
transformation $AW$ as the stochastic process with coordinates
$(AW)_h=W_{A^T h}$, for $h\in H_1$. If the process $W$ arises as
$W_h=\langle W,h\rangle_2$ from a random element $W$ in the Hilbert
space $H_2$, then this definition is consistent with identifying the
random element $A W$ in $H_1$ with the process $(\langle
AW,h\rangle_1\dvtx h\in H_1)$, as in (\ref{EqProcessG}) with $G=AW$.
Furthermore, if~$A$ is a \textit{Hilbert--Schmidt} operator (i.e., $AA^T$
is of trace class), and $W=Z$ is the iso-Gaussian process, then the
process $AW$ can be realized as a random variable in $H_1$ with a
$N(0,AA^T)$-distribution.

In the Bayesian setup the prior, which we take
$N(0,\Lambda)$, is the marginal distribution of $\mu$,
and the noise $Z$ in (\ref{EqProblem}) is
considered independent of $\m$. The joint distribution of $(Y,\m)$
is then also Gaussian, and so is the conditional distribution of
$\m$ given $Y$, the posterior distribution of $\m$.
In general, one must be a bit careful with manipulating possibly
``improper'' Gaussian distributions (see \cite{Mandelbaum}), but
in our situation the posterior is a
proper Gaussian conditional distribution on $H_1$.
%
%pr3.1 #&#
%
\begin{proposition}[(Full posterior)]\label{Posterior}
If $\m$ is $N(0,\Lambda)$-distributed and
$Y$ given $\m$ is $N(K\m, n^{-1}I)$-distributed, then
the conditional distribution of $\m$ given $Y$ is
Gaussian $N(AY, S_n)$ on $H_1$, where
%
%e3.2 ###
%e3.2 #&#
%
\begin{equation}\label{PostCov}
S_n = \Lambda-A(n^{-1}I+K\Lambda K^T)A^T,
\end{equation}
and $A\dvtx H_2\to H_2$ is the continuous linear operator
%
%e3.3 ###
%e3.3 #&#
%
\begin{equation}\label{A}
A=\Lambda^{1/2}\biggl(\frac1{n}I+\Lambda^{1/2}K^TK\Lambda
^{1/2}\biggr)^{-1}
\Lambda^{1/2}K^T=\Lambda K^T\biggl(\frac1nI+K\Lambda
K^T\biggr)^{-1}.\hspace*{-28pt}
\end{equation}
The posterior distribution is proper (i.e., $S_n$ has finite trace)
and equivalent (in the sense of absolute continuity) to the prior.
\end{proposition}
\begin{pf}
Identity (\ref{A}) is a special case of the
identity $(I+BB^T)^{-1} B=B(I+B^TB)^{-1}$, which is valid
for any compact, linear operator $B\dvtx H_1\to H_2$.
That $S_n$ is of trace class is a consequence of the fact
that it is bounded above by $\Lambda$ (i.e., $\Lambda-S_n$ is nonnegative
definite), which is of trace class by assumption.

The operator $\Lambda^{1/2}K^TK\Lambda^{1/2}\dvtx H_1\to H_1$
has trace bounded by $\|K^TK\|\tr(\Lambda)$ and hence
is of trace class. It follows that the variable
$\Lambda^{1/2}K^TZ$ can be defined as a random element in the
Hilbert space $H_1$, and so can $AY$, for~$A$ given by the first
expression in (\ref{A}). The joint distribution of $(Y,\m)$
is Gaussian with zero mean and covariance operator
\[
\pmatrix{
n^{-1}I + K\Lambda K^T & K \Lambda\cr
\Lambda K^T& \Lambda}.
\]
Using this with the second form of $A$ in (\ref{A}), we can
check that the cross covariance operator
of the variables $\m-AY$ and $Y$ (the latter viewed
as a~Gaussian stochastic process in $\RR^{H_2}$) vanishes and, hence,
these variables are independent.
Thus, the two terms in the decomposition $\m=(\m-AY) +AY$
are conditionally independent and degenerate given $Y$,
respectively. The distribution of $\m-AY$ is zero-mean
Gaussian with covariance operator $\Cov(\m-AY)=\Cov(\m)-\Cov(AY)$,
by the independence of $\m-AY$ and~$AY$.
This gives the form of the posterior distribution.

The final assertion may be proved by explicitly comparing
the Gaussian prior and posterior. Easier is to note that it
suffices to show that the model consisting of all
$N(K\m,n^{-1}I)$-distributions is dominated. In that case
the posterior can be obtained using Bayes' rule, which
reveals the normalized likelihood as a density relative to the
(in fact, \textit{any}) prior. To prove domination, we may consider
equivalently the distributions $\bigotimes_{i=1}^\infty N(\k_i\m_i,n^{-1})$
on $\RR^\infty$ of the sufficient statistic $(Y_i)$ defined
as the coordinates of $Y$ relative to the conjugate spectral
basis. These distributions, for $(\m_i)\in\ell_2$, are
equivalent to the distribution $\bigotimes_{i=1}^\infty N(0,n^{-1})$,
as can be seen with the help of Kakutani's theorem,
the affinity being $\exp(-\sum_i\k_i^2\m_i^2/8)>0$. (This
argument actually proves the well-known fact that
the Gaussian shift experiment obtained by translating
the standard normal distribution on $\RR^\infty$ over its RKHS
$\ell_2$ is dominated.)
\end{pf}

In the remainder of the paper we study the asymptotic behavior of the
posterior distribution, under the assumption that $Y= K\m_0 +
n^{-1/2}Z$ for a~fixed $\m_0 \in H_1$. The posterior is characterized
by its \textit{center} $AY$, the posterior mean, and its \textit{spread},
the posterior covariance operator $S_n$. The first depends on
the data, but the second is deterministic. From a frequentist-Bayes
perspective both are important: one would like the posterior mean
to give a good estimate for $\m_0$, and the spread to give a
good indication of the uncertainty in this estimate.

The posterior mean is a regularization,
of the Tikhonov type, of the naive estimator $K^{-1}Y$. It can also be
characterized as a penalized least squares estimator (see
\cite{Tikhonov,Mathe}): it
minimizes the functional
\[
\m\mapsto\|Y-K\m\|_2^2 + \frac1n\|\Lambda^{-1/2}\m\|^2_1.
\]
The penalty $\|\Lambda^{-1/2}\m\|_1$ is interpreted as $\infty$ if
$\m$ is not in the range of $\Lambda^{1/2}$. Because this range is
precisely\vspace*{1pt} the \textit{reproducing kernel Hilbert space} (RKHS) of the
prior (cf. \cite{vdVvZRKHS}),
with $\|\Lambda^{-1/2}\m\|_1$ as the RKHS-norm of $\m$, the
posterior mean also fits into the general regularization framework
using RKHS-norms (see \cite{Wahba}). In any case the posterior mean is a
well-studied point estimator in the literature on inverse problems.
In this paper we add a Bayesian interpretation to it, and are (more)
concerned with the full posterior distribution.

Next consider the posterior distribution of a linear functional $L\m$
of the parameter. We are not only interested in continuous, linear functionals
$L\m=\langle\m,l\rangle_1$, for some given $l\in H_1$, but also in
certain discontinuous functionals, such as point evaluation in a
Hilbert space of functions. The latter entail some technicalities. We
consider \textit{measurable linear functionals relative to the prior}
$N(0,\Lambda)$, defined in \cite{Skorohod}, pages 27--29, as Borel
measurable maps
$L\dvtx H_1\to\RR$ that are linear on a measurable linear subspace
$\underline H_1\subset H_1$ such that $N(0,\Lambda)(\underline
H_1)=1$. This definition is exactly right to make the marginal
posterior Gaussian.
%
%pr3.2 #&#
%
\begin{proposition}[(Marginal posterior)]\label{PosteriorLinear}
If $\m$ is $N(0,\Lambda)$-distributed and~%
$Y$ given $\m$ is $N(K\m, n^{-1}I)$-distributed, then
the conditional distribution\vadjust{\goodbreak} of $L\m$ given $Y$ for
a $N(0,\Lambda)$-measurable linear functional $L\dvtx H_1\to\RR$ is
a Gaussian distribution $N(LAY, s_n^2)$ on $\RR$, where
%
%e3.4 ###
%e3.4 #&#
%
\begin{equation}\label{PostCovL}
s_n^2 = (L\Lambda^{1/2})(L\Lambda^{1/2})^T-LA(n^{-1}I+K\Lambda K^T)(LA)^T,
\end{equation}
and $A\dvtx H_2\to H_2$ is the continuous linear operator defined
in (\ref{A}).
\end{proposition}
\begin{pf}
As in the proof of Proposition \ref{Posterior}, the first term
in the decomposition $L\m=L(\m-AY)+LAY$ is independent of $Y$.
Therefore, the posterior distribution is the marginal
distribution of $L(\m-AY)$ shifted by $LAY$. It suffices
to show that this marginal distribution is $N(0,s_n^2)$.

By Theorem 1 on page 28 in \cite{Skorohod}, there exists a sequence of
continuous linear maps $L_m\dvtx H_1\to\RR$ such that $L_mh\ra Lh$ for all
$h$ in a set with probability one under the prior
$\Pi=N(0,\Lambda)$. This implies that $L_m \Lambda^{1/2}h\ra
L\Lambda^{1/2}h$ for \textit{every} $h\in H_1$. Indeed, if $V=\{h\in
H_1\dvtx L_mh\ra
Lh\}$ and $g\notin V$, then $V_1:=V+ g$ and $V$ are disjoint measurable,
affine subspaces of $H_1$, where $\Pi(V)=1$. The range of $\Lambda^{1/2}$
is the RKHS of $\Pi$ and, hence, if $g$ is in this range, then
$\Pi(V_1)>0$, as $\Pi$ shifted over an element from its
RKHS is equivalent to $\Pi$. But then $V$ and $V_1$ are not disjoint.

Therefore, from the first definition of $A$ in (\ref{A}) we see that
$L_mA\ra LA$, and, hence, $L_m(\m-AY)\ra L(\m-AY)$, almost surely. As
$L_m$ is continuous, the variable $L_m(\m-AY)$ is normally distributed
with mean zero and variance
$L_mS_mL_m^T=(L_m\Lambda^{1/2})(L_m\Lambda^{1/2})^T-
L_mA(n^{-1}I+K\Lambda K^T)(L_mA)^T$, for $S_n$ given by
(\ref{PostCov}). The desired result follows upon taking the limit as
$m\ra\infty$.
\end{pf}

As shown in the preceding proof, $N(0,\Lambda)$-measurable linear
functionals $L$ automatically have the further property that
$L\Lambda^{1/2}\dvtx H_1\to\RR$ is a continuous linear map. This shows
that $LA$ and the adjoint operators $(L\Lambda^{1/2})^T$ and $(LA)^T$
are well defined, so that the formula for $s_n^2$ makes sense. If $L$
is a~\textit{continuous} linear operator, one can also write these
adjoints in terms of the adjoint\vspace*{1pt} $L^T$ of $L$, and express $s_n^2$ in
the covariance operator $S_n$ of Proposition~\ref{Posterior} as
$s_n^2=LS_nL^T$. This is exactly as expected.

In the remainder of the paper we study the full posterior distribution
$N(AY,S_n)$, and its marginals $N(LAY,s_n^2)$. We are particularly
interested in the influence of the prior on the performance of the
posterior distribution for various true parameters~$\m_0$. We study
this in the following setting.
%
%as3.1 #&#
%
\begin{assumption}\label{common}\label{PPP}
The operators $K^T K$ and $\Lambda$ have the same
eigenfunctions~$(e_i)$, with eigenvalues $(\k_i^2)$ and $(\lambda_i)$,
satisfying
%
%e3.5 ###
%e3.5 #&#
%
\begin{equation}\label{SetupPPP}
\lambda_i = \t_n^2i^{-1-2\a},\qquad C^{-1}i^{-p} \leq\k_i \leq C i^{-p}
\end{equation}
for some $\a>0 $, $p \geq0$, $C\ge1$ and $\t_n > 0$ such that
$n\t_n^2\ra\infty$.
Furthermore, the true parameter $\m_0$ belongs to $S^\b$ for some $\b>0$:
that is, its coordinates~$(\m_{0,i})$ relative to $(e_i)$ satisfy
$\sum_{i=1}^\infty\m_{0,i}^2i^{2\b}<\infty$.
\end{assumption}

The setting of Assumption \ref{PPP} is a Bayesian extension of the
\textit{mildly ill-posed inverse problem} (cf. \cite{Cavalier}). We
refer to the parameter $\b$ as the ``regularity'' of the
true parameter $\m_0$. In the special case that $H_1$ is a function
space and $(e_i)$ its Fourier basis, this parameter gives smoothness
of $\m_0$ in the classical Sobolev sense. Because the coefficients
$(\m_i)$ of the prior parameter $\m$ are normally
$N(0,\l_i)$-distributed, under Assumption \ref{PPP} we have $\E
\sum_ii^{2\a'}\m_i^2=\t_n^2\sum_ii^{2\a'}\l_i<\infty$ if and
only if
$\a'<\a$. Thus, $\a$ is ``almost'' the smoothness of the parameters
generated by the prior. This smoothness is modified by the
\textit{scaling factor} $\t_n$. Although this leaves the relative sizes
of the coefficients $\m_i$, and hence the qualitative smoothness of
the prior, invariant, we shall see that scaling can completely alter
the performance of the Bayesian procedure. Rates $\t_n\da0$ increase,
and rates $\t_n\ua\infty$ decrease the regularity.

%s4 ###
%se4 #&#
\section{Recovering the full parameter}
\label{SectionFull}
We denote by $\Pi_n(\cdot\,\given Y)$ the posterior distribution
$N(AY,S_n)$, derived in Proposition \ref{Posterior}. Our first
theorem shows that it contracts as $n\ra\infty$ to the true parameter
at a rate $\eps_n$ that depends on all four parameters $\a,\b,\t_n,p$
of the (Bayesian) inverse problem.
%
%th4.1 #&#
%
\begin{theorem}[(Contraction)]\label{ContrPPP}
If $\m_{0}$, $(\lambda_i)$, $(\k_i)$ and $(\t_n)$ are
as in Assumption~\ref{PPP}, then
$\E_{\m_0}\Pi_n(\m\dvtx \|\m-\m_0\|_1 \geq M_n \eps_n \given
Y) \ra0$,
for every $M_n\ra\infty$, where
%
%e4.1 ###
%e4.1 #&#
%
\begin{equation}\label{Rate}
\eps_n =(n\t_n^2)^{-{\b}/({1+2\a+2p})\wedge1} + \t_n(n\t
_n^2)^{-{\a}/({1+2\a+2p})}.
\end{equation}
The rate is uniform over $\m_0$ in balls in $S^\b$. In particular:
\begin{longlist}[(iii)]
\item[(i)] If $\t_n \equiv1$, then $\eps_n = n^{-{(\a\wedge\b
)}/{(1+2\a+2p)}}$.
\item[(ii)] If $\b\leq1+2\a+ 2p$ and $\t_n \asymp n^{(\a-\b
)/(1+2\b+2p)}$,
then $\eps_n = n^{-\b/(1+2\b+2p)}$.
\item[(iii)] If $\b> 1+2\a+ 2p$, then $\eps_n \gg n^{-\b/(1+2\b+2p)}$,
for every scaling $\t_n$.
\end{longlist}
\end{theorem}

The minimax rate of convergence over a Sobolev ball $S^\b$ is of the
order $n^{-\b/(1+2\b+2p)}$ (see \cite{Cavalier}). By (i) of the
theorem the
posterior contraction rate is the same
if the regularity of the prior is chosen to match the
regularity of the truth ($\a=\b$) and the scale $\t_n$ is fixed.
Alternatively, the optimal rate is also attained by appropriately
scaling ($\t_n \asymp n^{(\a-\b)/(1+2\b+2p)}$, determined
by balancing the two terms in ${\varepsilon}_n$) a prior that is
regular enough ($\b\le1+2\a+2p$). In all other cases
(no scaling and $\a\not=\b$, or any scaling combined with a
rough prior $\b> 1+2\a+ 2p$), the contraction rate is
slower than the minimax rate.

That ``correct'' specification of the prior gives the optimal rate is
comforting to the true Bayesian. Perhaps the main message of the
theorem is that even if the prior mismatches the truth, it may be
scalable to give the optimal rate. Here, similar as found by
\cite{vdVvZRescaledGaussian} in a different setting, a smooth prior
can be scaled to make it ``rougher'' to any degree, but a rough prior
can be ``smoothed'' relatively little (namely, from $\a$ to any
$\b\le1+2\a+2p$). It will be of interest to investigate a full or
empirical Bayesian approach to set the scaling parameter.

Bayesian inference takes the spread in the posterior distribution as
an expression of uncertainty. This practice is not validated by (fast)
contraction of the posterior. Instead we consider the
frequentist coverage of credible sets. As the posterior
distribution is Gaussian, it is natural to center a~credible region
at the posterior mean. Different shapes of such a set could be
considered. The natural counterpart of the preceding theorem is to consider
balls. In the next section we also consider bands. (Alternatively, one
might consider ellipsoids, depending on geometry of the support of the
posterior.)

Because the posterior spread $S_n$ is deterministic, the radius is the
only degree of freedom when we choose a ball, and we fix it by the
desired ``credibility level'' $1-\g\in(0,1)$. A \textit{credible
ball} centered at the posterior mean $AY$ takes the form, where
$B(r)$ denotes a ball of radius $r$ around 0,
%
%e4.2 ###
%e4.2 #&#
%
\begin{equation}
\label{EqCredReg}
AY+B(r_{n,\g}):=\{\m\in H_1 \dvtx \|\m-AY\|_1<r_{n,\g}\},
\end{equation}
where the radius $r_{n,\g}$ is determined so that
%
%e4.3 ###
%e4.3 #&#
%
\begin{equation}
\label{EqRadius}
\Pi_n\bigl(AY+B(r_{n,\g})\given Y\bigr)=1-\g.
\end{equation}
Because the posterior spread $S_n$ is not dependent on the data,
neither is the radius~$r_{n,\g}$.
The frequentist \textit{coverage} or confidence of the set
(\ref{EqCredReg}) is
%
%e4.4 ###
%e4.4 #&#
%
\begin{equation}
\label{EqCoverage}
\Pr_{\m_0}\bigl( \m_0\in AY+B(r_{n,\g})\bigr),
\end{equation}
where under the probability measure $\Pr_{\m_0}$
the variable $Y$ follows (\ref{EqProblem}) with $\m=\m_0$. We shall consider
the coverage as $n\ra\infty$ for fixed $\m_0$, uniformly
in Sobolev balls, and also along sequences $\m_0^n$ that change
with $n$.

The following theorem shows that the relation of the coverage to the
credibility level $1-\g$ is mediated by all parameters of the
problem. For further insight, the credible region is also compared
to the ``correct'' frequentist confidence ball $AY+B(\wt r_{n,\g})$,
which has radius $\wt r_{n,\g}$ chosen so that the probability in
(\ref{EqCoverage}) with $r_{n,\g}$ replaced by $\wt r_{n,\g}$ is
equal to $1-\g$.
%
%th4.2 #&#
%
\begin{theorem}[(Credibility)]\label{CrS}\label{CrNS}
Let $\m_{0}$, $(\l_i)$, $(\k_i)$, and $\t_n$
be as in Assumption~\ref{PPP}, and set $\wt\b= \b\wedge(1+2\a+2p)$.
The asymptotic coverage of the credible region (\ref{EqCredReg}) is:
\begin{longlist}[(iii)]
\item[(i)] 1, uniformly in $\m_0$ with $\|\m_0\|_\b\leq1$,
if $\t_n \gg n^{(\a-\wt\b)/(1+2\wt\b+2p)}$;
in this case $\wt r_{n,\g}\asymp r_{n,\g}$.
\item[(ii)] 1, for every fixed $\m_0 \in S^\b$, if $\b< 1+2\a+2p$ and
$\t_n \asymp n^{(\a-\wt\b)/(1+2\wt\b+2p)}$;
$c$,
along some $\m_{0}^{n}$ with ${\sup_n }\|\m_0^{n}\|_\b< \infty$,
if $\t_n \asymp n^{(\a-\wt\b)/(1+2\wt\b+2p)}$ (any $c\in[0,1)$).
\item[(iii)] 0, along some $\m_0^n$ with ${\sup_n }\|\m_0^{n}\|_\b<
\infty$,\vadjust{\goodbreak}
if $\t_n \ll n^{(\a-\wt\b)/(1+2\wt\b+2p)}$.
\end{longlist}
If $\t_n\equiv1$, then the cases \textup{(i), (ii)} and \textup{(iii)} arise if
$\a<\b$, $\a=\b$ and $\a>\b$, respectively. In case \textup{(iii)} the sequence
$\m_0^n$ can then be chosen fixed.
\end{theorem}

The theorem is easiest to interpret in the situation without scaling
(\mbox{$\t_n\equiv1$}). Then oversmoothing the prior [case (iii): $\a>\b$] has
disastrous consequences for the coverage of the credible sets,
whereas undersmoothing [case~(i): $\a<\b$] leads to conservative
confidence sets. Choosing a prior of correct regularity [case (ii):
$\a=\b$] gives mixed results.

Inspection of the proofs shows that the lack of coverage in case of
oversmoothing arises from a bias in the positioning of the posterior
mean combined with a posterior spread that is smaller even than in the
optimal case. In other words, the posterior is off mark, but believes
it is very right. The message is that (too) smooth priors
should be avoided; they lead to overconfident posteriors, which
reflect the prior information rather than the data, even if the amount
of information in the data increases indefinitely.

Under- and correct smoothing give very
conservative confidence regions (coverage equal to 1). However,
(i) and (ii) also
show that the credible ball has the same order of magnitude
as a correct confidence ball ($1\ge\wt r_{n,\g}/\allowbreak r_{n,\g}\gg0$), so
that the spread in the posterior does give the correct \textit{order} of
uncertainty. This at first sight surprising phenomenon is caused
by the fact that the posterior distribution concentrates near the boundary
of a~ball around its mean, and is not spread over the inside of the ball.
The coverage is 1, because this sphere is larger than
the corresponding sphere of the frequentist distribution of $AY$, even though
the two radii are of the same order.

By Theorem \ref{ContrPPP} the optimal contraction rate is obtained
(only) by a prior of the correct smoothness. Combining the two
theorems leads to the conclusion that priors that slightly undersmooth
the truth might be preferable. They attain a nearly optimal rate of contraction
and the spread of their posterior gives a reasonable sense of uncertainty.

Scaling\vspace*{1pt} of the prior modifies these conclusions. The optimal scaling
$\t_n\asymp n^{(\a-\b)/(1+2\a+2p)}$ found in Theorem \ref{ContrPPP},
possible if $\b<1+2\a+2p$, is covered in case (ii). This rescaling
leads to a balancing of square bias, variance and spread, and to
credible regions of the correct order of magnitude, although the
precise (uniform) coverage can be any number in $[0,1)$.
Alternatively, bigger rescaling rates are covered in case (i)
and lead to coverage 1. The optimal or slightly bigger rescaling
rate seems the most sensible. It would be interesting to extend
these results to data-dependent scaling.

%s5 ###
%se5 #&#
\section{Recovering linear functionals of the parameter}
\label{SectionLinear}
We denote by $\Pi_n(\m\dvtx\break L\m\in\cdot\,\given Y)$ the posterior
distribution of the linear functional $L$, as described in
Proposition \ref{PosteriorLinear}. A continuous, linear functional
$L\dvtx
H_1\to\RR$ can be identified with an inner product $L\m=\langle\m,
l\rangle_1$, for some $l\in H_1$, and hence with a~sequence $(l_i)$ in
$\ell_2$ giving its coordinates in the eigenbasis $(e_i)$.

As shown in the proof of Proposition \ref{PosteriorLinear}, for $L$ in
the larger class of $N(0,\Lambda)$-measurable linear functionals, the
functional $L\Lambda^{1/2}$ is a continuous linear map on $H_1$
and hence can be identified with an element of
$H_1$. For such a functional $L$ we define a sequence $(l_i)$ by
$l_i=(L\Lambda^{1/2})_i/\sqrt{\l_i}$, for
$((L\Lambda^{1/2})_i)$ the coordinates of $L\Lambda^{1/2}$
in the eigenbasis. The assumption that~$L$ is a~%
$N(0,\Lambda)$-measurable linear functional implies that
$\sum_il_i^2\l_i<\infty$, but $(l_i)$ need not be contained
in $\ell_2$; if \mbox{$(l_i)\in\ell_2$}, then $L$ is
continuous and the definition of $(l_i)$ agrees with the definition in
the preceding paragraph.

We measure the smoothness of the functional $L$ by the size
of the coefficients~$l_i$, as $i\ra\infty$. First we assume that
the sequence is in $S^q$, for some~$q$.

%th5.1 #&#
%
\begin{theorem}[(Contraction)]
\label{LinContrPPP}
If $\m_{0}$, $(\l_i)$, $(\k_i)$ and $(\t_n)$ are
as in Assumption~\ref{PPP} and the representer $(l_i)$ of
the $N(0,\Lambda)$-measurable linear functional $L$ is contained in
$S^q$ for
$q\ge-\b$, then
$\E_{\m_0}\Pi_n(\m\dvtx |L\m- L\m_0| \ge\break M_n \eps_n \given
Y) \ra0$,
for every sequence $M_n\to\infty$, where
\[
\eps_n = (n\t_n^2)^{-({\b+q})/({1+2\a+2p})\wedge1}
+\t_n(n\t_n^2)^{-({1/2+\a+q})/({1+2\a+2p})\wedge(1/2)}.
\]
The rate is uniform over $\m_0$ in balls in $S^\b$. In particular:
\begin{longlist}[(iii)]
\item[(i)] If $\t_n\equiv1$, then
$\eps_n=n^{-(\b\wedge(1/2+\a)+q)/(1+2\a+2p)}\vee n^{-1/2}$.
\item[(ii)] If\vspace*{1pt} $q \leq p$ and $\b+q \leq1 + 2\a+ 2p$ and
$\t_n \asymp n^{({1/2+\a-\b})/({2\b+ 2p})}$, then
$\eps_n = n^{-(\b+q)/(2\b+2p)}$.
\item[(iii)] If $q\le p$ and $\b+ q > 1+2\a+2p$, then
$\eps_n \gg n^{-({\b+q})/({2\b+2p})}$ for every scaling $\t_n$.
\item[(iv)] If $q \geq p$ and
$\t_n\gtrsim n^{({1/2+\a-\wt\b+ p-q})/({2\wt\b+ 2q})}$,
where $\wt\b= \b\wedge(1+2\a+2p-q)$, then $\eps_n = n^{-1/2}$.
\end{longlist}
\end{theorem}

If $q\ge p$, then the smoothness of the functional $L$ cancels the
ill-posedness of the operator $K$, and estimating $L\m$ becomes a
``regular'' problem with an~$n^{-1/2}$ rate of convergence. Without
scaling the prior ($\tau_n\equiv1$), the posterior contracts at this
rate [see (i) or (iv)] if the prior is not too smooth
$(\a\le\b-1/2+q-p$). With scaling, the rate is also attained, with any
prior, provided the scaling parameter $\t_n$ does not tend to zero too
fast [see (iv)]. Inspection of the proof shows that too smooth priors
or too small scale creates a bias that slows the rate.

If $q<p$, where we take $q$ the ``biggest'' value such that $(l_i)\in
S^q$, estimating~$L\m$ is still an inverse problem. The minimax rate
over a ball in the Sobolev space $S^\b$ is known to be bounded above
by $n^{-(\b+q)/(2\b+2p)}$ (see
\cite{Donoho,DonohoLow,Goldenshluger}).

This rate is attained without scaling [see (i): $\tau_n\equiv1$] if
and only if the prior smoothness $\a$ is equal to the true smoothness
$\b$ minus $1/2$ ($\a+1/2=\b$). An intuitive explanation for this
apparent mismatch of prior and truth is that regularity of the
parameter in the Sobolev scale ($\m_0\in S^\b$) is not the appropriate
type of regularity for estimating a linear functional $L\m$. For
instance, the difficulty of estimating a function at a point is
determined by the regularity in a neighborhood of the point, whereas
the Sobolev scale measures global regularity over the domain. The fact
that a Sobolev space of order $\b$ embeds continuously in a H\"older
space of regularity $\b-1/2$ might give a quantitative explanation of
the ``loss'' in smoothness by $1/2$ in the special case that the
eigenbasis is the Fourier basis.
In our Bayesian context we draw the conclusion that the prior must be
adapted to the inference problem if we want to obtain the optimal
frequentist rate: for estimating the global parameter, a
good prior must match the truth ($\a=\b$), but for estimating a linear
functional a good prior must consider a Sobolev truth of order $\b$ as
having regularity $\a=\b-1/2$.

If the prior smoothness $\a$ is not $\b-1/2$, then the minimax rate
may still be attained by scaling the prior. As in the global problem,
this is possible only if the prior is not too rough [$\b+q\le1+2\a+2p$,
cases (ii) and (iii)]. The optimal scaling when this is possible [case
(ii)] is the same as the optimal scaling for the global problem
[Theorem \ref{ContrPPP}(ii)] \textit{after} decreasing $\b$ by $1/2$.
So the
``loss in regularity'' persists in the scaling rate. Heuristically
this seems to imply that a simple data-based procedure to set the
scaling, such as empirical or hierarchical Bayes, cannot attain
simultaneous optimality in both the global and local senses.

In the application of the preceding theorem, the functional $L$, and
hence the sequence $(l_i)$, will be given. Naturally, we apply the
theorem with $q$ equal to the largest value such that $(l_i)\in
S^q$. Unfortunately, this lacks precision for the sequences $(l_i)$
that decrease exactly at some polynomial order: a~sequence $l_i\asymp
i^{-q-1/2}$ is in $S^{q'}$ for every $q'<q$, but not in $S^q$. In the
following theorem we consider these sequences, and the slightly more
general ones such that $|l_i|\le i^{-q-1/2}\S(i)$, for some slowly
varying sequence $\S(i)$. Recall that $\S\dvtx [0,\infty)\to\RR$ is
\textit{slowly varying} if $\S(tx)/\S(t)\ra1$ as $t\ra\infty$, for
every $x>0$. [For these sequences $(l_i)\in S^{q'}$ for every
$q'<q$, $(l_i)\notin S^{q'}$ for $q'>q$,
and $(l_i)\in S^q$ if and only if $\sum_i\S^2(i)/i<\infty$.]
%
%th5.2 #&#
%
\begin{theorem}[(Contraction)]\label{LinContrPPPRV}
If $\m_{0}$, $(\l_i)$, $(\k_i)$ and $(\t_n)$ are
as in Assumption~\ref{PPP} and the representer $(l_i)$ of
the $N(0,\Lambda)$-measurable linear functional $L$
satisfies $|l_i|\le i^{-q-1/2}\S(i)$ for a slowly varying function $\S
$ and
$q>-\b$, then the result of Theorem \ref{LinContrPPP} is valid with
%
%e5.1 ###
%e5.1 #&#
%
\begin{equation}
\label{EqTauLRV}
\eps_n = (n\t_n^2)^{-({\b+q})/({1+2\a+2p})\wedge1}\g_n
+\t_n(n\t_n^2)^{-({1/2+\a+q})/({1+2\a+2p})\wedge(1/2)}\d_n,\hspace*{-35pt}
\end{equation}
where, for $\r_n=(n\t_n^2)^{1/(1+2\a+2p)}$,
\begin{eqnarray*}
\g_n^2&=&\cases{ \S^2(\r_n), &\quad if $\b+q < 1+2\a+2p$,\vspace*{2pt}\cr
\displaystyle \sum_{i \leq\r_n}\frac{\S^2(i)}{i}, &\quad if $\b+q=1+2\a+2p$,\vspace*{2pt}\cr
1, &\quad if $\b+q>1+2\a+2p$,}\\
\d_n^2&=&\cases{\S^2(\r_n), &\quad if $q < p$,\vspace*{2pt}\cr
\displaystyle \sum_{i \leq\r_n}\frac{\S^2(i)}{i}, &\quad if $q=p$,\vspace*{2pt}\cr
1, &\quad if $q>p$.}
\end{eqnarray*}
This has the same consequences as in Theorem \ref{LinContrPPP},
up to the addition of slowly varying terms.
\end{theorem}

Because the posterior distribution for the linear
functional $L\m$ is the one-dimensional normal distribution $N(LAY,
s_n^2)$, the natural \textit{credible interval} for $L\m$ has endpoints
$LAY \pm z_{\g/2}s_n$, for $z_{\g}$ the (lower)
standard normal $\g$-quantile. The
\textit{coverage} of this interval is
\[
\Pr_{\m_0} (LAY+z_{\g/2}s_n\le L\m_0 \leq LAY -z_{\g
/2}s_n),
\]
where $Y$ follows (\ref{EqProblem}) with $\m=\m_0$.
To obtain precise results concerning coverage, we assume that $(l_i)$
behaves polynomially up to a slowly varying term, first in the
situation $q<p$ that estimating $L\m$ is an inverse problem. Let $\wt
\t_n$ be the (optimal) scaling
$\t_n$ that equates the two terms in the right-hand side of (\ref{EqTauLRV}).
This satisfies $\wt\t_n\asymp n^{({1/2+\a-\wt\b
})/({2\wt\b+ 2p})}\h_n$, for a slowly varying factor~$\h_n$, where
$\wt\b=\b\wedge(1+2\a+2p-q)$.
%
%th5.3 #&#
%
\begin{theorem}[(Credibility)]
\label{LinCrS}
Let $\m_{0}$, $(\l_i)$, $(\k_i)$ and $(\t_n)$ be as in
Assumption~\ref{PPP}, and let $|l_i|=i^{-q-1/2}\S(i)$
for $q < p$ and a slowly varying function~$\S$. Then the
asymptotic coverage of the interval $LAY \pm z_{\g/2}s_n$ is:
\begin{longlist}[(iii)]
\item[(i)] in $(1-\g,1)$, uniformly in $\m_0$ such that $\|\m_0\|
_\b\leq1$
if $\t_n \gg\wt\t_n$.
\item[(ii)] in $(1-\g,1)$, for every $\m_0\in S^\b$, if $\t
_n\asymp\wt\t_n$
and $\b+q < 1+2\a+2p$;
in $(0,c)$, along some $\m_0^{n}$ with $\sup_n\|\m_0^{n}\|
_\b< \infty$
if $\t_n\asymp\wt\t_n$ [any $c\in(0, 1)$].
\item[(iii)] $0$ along some $\m_0^n$ with $\sup_n\|\m_0^{n}\|_\b<
\infty$
if $\t_n \ll\wt\t_n$.
\end{longlist}
In case \textup{(iii)} the sequence $\m_0^n$ can be taken a fixed element
$\m_0$ in $S^\b$ if $\t_n\lesssim n^{-\d}\wt\t_n$ for some $\d>0$.

Furthermore, if $\t_n\equiv1$, then the coverage takes the form
as in \textup{(i), (ii)} and~\textup{(iii)} if $\a<\b-1/2$, $\a=\b-1/2$, and $\a>\b-1/2$,
respectively, where in case~\textup{(iii)} the sequence $\m_0^n$ can be taken
a fixed element.
\end{theorem}

Similarly, as in the nonparametric problem, oversmoothing leads to
coverage 0, while undersmoothing gives conservative intervals.
Without scaling the cut-off for under- or oversmoothing is at
$\a=\b-1/2$; with scaling the cut-off for the scaling rate is at the
optimal rate $\wt\t_n$.

The conservativeness in the case of undersmoothing is less
extreme for functionals than for the full parameter, as the coverage
is strictly between the credibility level $1-\g$ and 1.
The general message is the same: oversmoothing is disastrous for
the interpretation of credible band, whereas undersmoothing
gives bands that at least have the correct order of magnitude, in the
sense that their width is of the same order as the variance of the
posterior mean (see the proof). Too much undersmoothing is also
undesirable, as it leads to very wide confidence bands, and may cause
that $\sum_i l_i^2\lambda_i$ is no longer finite (see measurability property).

The results (i) and (ii) are the same for every $q < p$, even if $\tau
_n\equiv1$. Closer inspection would reveal that for a given $\mu_0$
the exact coverage depends on $q$ [and $\S(i)$] in a complicated way.

If $q\ge p$, then the smoothness of the functional $L$ compensates the
lack of smoothness of $K^{-1}$, and estimating $L\m$ is not a true
inverse problem. This drastically changes the performance of
credible intervals. Although oversmoothing again destroys their
coverage, credible intervals are exact confidence sets if the prior
is not too smooth. We formulate this in terms of a~Bernstein--von Mises
theorem.

The Bernstein--von Mises theorem for parametric models asserts that the
posterior distribution approaches a normal distribution centered at an
efficient estimator of the parameter and with variance equal to its
asymptotic variance. It is the ultimate link between Bayesian and
frequentist procedures. There is no version of this theorem
for infinite-dimensional parameters \cite{Freedman}, but the
theorem may hold for ``smooth'' finite-dimensional projections,
such as the linear functional $L\m$ (see~\cite{Castillo}).

In the present situation the posterior distribution of $L\m$ is
already normal by the normality of the model and the prior:
it is a $N(LAY,s_n^2)$-distribution by Proposition \ref{PosteriorLinear}.
To speak of a Bernstein--von Mises theorem, we also require the following:
\begin{longlist}[(iii)]
\item[(i)] That the (root of the) spread $s_n$ of the posterior distribution
is asymptotically equivalent to the standard deviation
$t_n$ of the centering variable $LAY$.
\item[(ii)] That the sequence $(LAY-L\m_0)/t_n$ tends in distribution to
a standard normal distribution.
\item[(iii)] That the centering $LAY$ is
an asymptotically efficient estimator of~$L\m$.
\end{longlist}
We shall show that (i) happens if and only if the
functional $L$ cancels the ill-posedness of the operator $K$,
that is, if $q\ge p$ in Theorem \ref{LinContrPPPRV}. Interestingly, the
rate of convergence $t_n$ must be $n^{-1/2}$ up to a slowly varying factor
in this case, but it could be strictly slower than $n^{-1/2}$ by a
slowly varying factor increasing to infinity.

Because $LAY$ is normally distributed by the normality of the model,
assertion (ii) is equivalent to saying that its bias tends to zero
faster than~$t_n$. This happens provided the prior does not oversmooth
the truth too much. For very smooth functionals ($q>p$) there is some
extra ``space'' in the cut-off for the smoothness, which (if the prior
is not scaled: $\t_n\equiv1$) is at $\a=\b-1/2+q-p$, rather than at
$\a=\b-1/2$ as for the (global) inverse estimating problem. Thus, the
prior may
be considerably smoother than the truth if the functional is very smooth.

Let \mbox{$\|\cdot\|$} denote the total variation norm between measures.
Say that $l\in\RV^q$ if $|l_i|=i^{-q-1/2}\S(i)$ for a slowly
varying function $\S$. Write
\[
B_n=\sup_{\|\m\|_\b\lesssim1}|LAK\m-L\m|
\]
for the maximal bias
of $LAY$ over a ball in the Sobolev space $S^\b$. Finally,
let $\wt\t_n$ be the (optimal) scaling $\t_n$
in that it equates the two terms in the right-hand side
of~(\ref{EqTauLRV}).
%
%th5.4 #&#
%
\begin{theorem}[(Bernstein--von Mises)]\label{TBvM-reg}
Let $\m_{0}$, $(\l_i)$, and $(\k_i)$ be as in Assumption \ref{PPP},
and let $l$ be the representer of the $N(0,\Lambda)$-measurable
linear functional~$L$:
\begin{longlist}[(iii)]
\item[(i)] If $l\in S^p$, then $s_n/t_n\ra1$; in this case
$nt_n^2\ra\sum_i l_i^2/\k_i^2$.
If $l\in\RV^q$, then $s_n/t_n\ra1$ if and only if $q\ge p$;
in this case $n\mapsto nt_n^2$ is slowly varying.
\item[(ii)] If $l\in S^q$ for $q\ge p$, then $B_n=o(t_n)$
if either $\t_n\gg n^{(\a+1/2-\b)/(2\b+2q)}$ or
($\t_n\equiv1$ and $\a<\b-1/2+q-p$).
If $l\in\RV^q$ for $q\ge p$, then $B_n=o(t_n)$ if
($\t_n \gg\wt\t_n$) or ($\t_n\equiv1$ and $\a<\b-1/2+q-p$)
or ($q=p$, $\t_n\equiv1$ and $\a=\b-1/2+q-p$)
or [$q> p$, $\t_n\equiv1$ and $\a=\b-1/2+q-p$ and $\S(i)\ra0$ as
$i\ra\infty$].
\item[(iii)] If $l\in S^p$ or $l\in\RV^p$ and $B_n=o(t_n)$,
then $\E_{\m_0}\|\Pi_n(L\m\in\cdot\,\given Y)-N(LAY$, $
t_n^2)\|\ra0$
and $(LAY-L\m_0)/t_n$ converges under $\m_0$
in distribution to a standard normal distribution, uniformly
in $\|\m_0\|_\b\lesssim1$. If $l\in S^p$, then this is also
true with $LAY$ and $t_n^2$ replaced by $\sum_iY_il_i/\k_i$ and
its variance $n^{-1}\sum_il_i^2/\k_i^2$.
\end{longlist}
In both cases \textup{(iii)},
the asymptotic coverage of the credible interval $LAY \pm z_{\g/2}s_n$
is $1-\g$, uniformly in $\|\m_0\|_\b\lesssim1$.
Finally, if the conditions under~\textup{(ii)} fail,
then there exists $\m_0^n$ with $\sup_n\|\m_0^n\|_\b<\infty$
along which the coverage tends to an arbitrarily low value.
\end{theorem}

The observation $Y$ in (\ref{EqProblem}) can be viewed as a reduction
by sufficiency of a random sample of size $n$ from the distribution
$N(K\m,I)$. Therefore, the model fits in the framework of
i.i.d. observations, and ``asymptotic efficiency'' can be defined in
the sense of semiparametric models discussed in, for example,
\cite{BKRW,vdV88} and \cite{vdVAS}. Because the model is shift-equivariant,
it suffices to consider local efficiency at $\m_0=0$. The
one-dimensional submodels $N(K(th),I)$ on the sample space
$\RR^{H_2}$, for $t\in\RR$ and a fixed ``direction'' $h\in H_1$, have
likelihood ratios
\[
\log\frac{dN(tKh,I)}{dN(0,I)}(Y)=tY_{Kh}-\frac12 t^2\|Kh\|_2^2.
\]
Thus, their \textit{score function} at $t=0$ is the $(Kh)$th coordinate
of a single observation $Y=(Y_h\dvtx h\in H_2)$, the \textit{score
operator} is the map $\tilde K\dvtx H_1\to L_2(N(0,I))$ given by
$\tilde Kh(Y)=Y_{Kh}$, and the \textit{tangent space} is the range of~%
$\tilde K$. [We denote the score operator by the same symbol~$K$ as in
(\ref{EqProblem}); if the observation $Y$ \textit{were} realizable
in~$H_2$, and not just in the bigger sample space~$\RR^{H_2}$, then
$Y_{Kh}$ would correspond to $\langle Y, Kh\rangle_2$ and, hence, the
score would be exactly $Kh$ for the operator in (\ref{EqProblem}) after
identifying $H_2$ and its dual space.] The adjoint of the score
operator restricted to the closure of the tangent space is the operator
$\tilde K^T\dvtx \overline{\tilde KH_1}\subset L_2 (N(0,I))\to H_1$
that satisfies $\tilde K^T(Y_g)=K^Tg$, where~$K^T$ on the right is the
adjoint of \mbox{$K\dvtx H_1\to H_2$}. The functional $L\m=\langle
l,\m\rangle_1$ has derivative $l$. Therefore, by \cite{vdV91}
asymptotically regular sequences of estimators exist, and the local
asymptotic minimax bound for estimating $L\m$ is finite, if and only if
$l$ is contained in the range of $K^T$. Furthermore, the variance bound
is $\|m\|_2^2$ for $m\in H_2$ such that $K^Tm=l$.

In our situation the range of $K^T$ is $S^p$, and if $l\in S^p$, then
by Theorem \ref{TBvM-reg}(iii) the variance of the posterior is
asymptotically equivalent to the variance bound and its centering can
be taken equal to the estimator $n^{-1}\sum_iY_il_i/\k_i$, which
attains this variance bound. Thus, the theorem gives a semiparametric
Bernstein--von Mises theorem, satisfying every of (i), (ii),~(iii) in
this case. If only $l\in\RV^p$ and not $l\in S^p$, the theorem still
gives a~Bernstein--von Mises type theorem, but the rate of convergence is
slower than $n^{-1/2}$, and the standard efficiency theory does not
apply.

%s6 ###
%se6 #&#
\section{Example---Volterra operator}
\label{SectionVolterra}
The classical \textit{Volterra operator} $K\colon L^2[0$, $1] \to L^2[0,1]$ and
its adjoint $K^T$ are given by
\[
K\m(x) = \int_0^x \m(s) \,ds,\qquad
K^T\m(x) = \int_x^1 \m(s) \,ds.
\]
The resulting problem (\ref{EqProblem}) can also be written in
``signal in white noise'' form as follows: observe the process
$(Y_t\dvtx
t\in[0,1])$
given by $Y_t=\break\int_0^t\int_0^s \m(u)\, du \,ds+n^{-1/2}W_t$, for a Brownian
motion $W$.

The eigenvalues, eigenfunctions of $K^TK$ and conjugate basis
are given by (see~\cite{Halmos}), for $i=1,2,\ldots,$
\begin{eqnarray*}
\k_i ^2&=& \frac{1}{(i-1/2)^2\pi^2},\qquad
e_i(x) = \sqrt{2}\cos\bigl((i-1/2)\pi x\bigr),\\
f_i(x) &=& \sqrt{2}\sin\bigl((i-1/2)\pi x\bigr).
\end{eqnarray*}
The $(f_i)$ are the eigenfunctions of $KK^T$,
relative to the same eigenvalues, and
$Ke_i = \k_i f_i$ and $K^T f_i = \k_ie_i$, for every $i\in\NN$.

To illustrate our results with simulated data, we start by
choosing a true function~$\m_0$, which we expand as $\m_0 = \sum_i
\m_{0,i}e_i$ on the basis $(e_i)$. The data are the function
\[
Y = K\m_0 + \frac1{\sqrt n}Z = \sum_i \m_{0,i}\k_i f_i + \frac
1{\sqrt n}Z.
\]
It can be generated relative to the conjugate basis $(f_i)$ as a
sequence of independent Gaussian random variables $Y_1, Y_2, \ldots$
with $Y_i \sim N(\m_{0,i}\k_i, n^{-1/2})$. The posterior
distribution of $\m$ is Gaussian with mean $AY$ and covariance
operator $S_n$, given in Proposition \ref{Posterior}. Under
Assumption \ref{PPP} it can be represented in terms of
the coordinates $(\m_i)$ of $\m$ relative to the basis $(e_i)$
as (conditionally) independent Gaussian variables $\m_1, \m_2, \ldots
$ with
\[
\m_i \big|Y \sim N\biggl(\frac{n\l_i\k_iY_i}{1+n\l_i\k_i^2},\frac{\l
_i}{1+n\l_i\k_i^2}\biggr).
\]
The (marginal) posterior distribution for the function $\m$ at a point
$x$ is obtained by expanding $\m(x) = \sum_i \m_ie_i(x)$, and applying
the framework of linear functionals $L\m=\sum_il_i\m_i$
with $l_i =e_i(x)$. This
shows that
\[
\m(x)\big|Y \sim N\biggl(\sum_i \frac{n\l_i\k_iY_ie_i(x)}{1+n\l_i\k
_i^2},\sum_i \frac{\l_ie_i(x)^2}{1+n\l_i\k_i^2}\biggr).
\]
We obtained (marginal) posterior credible bands by
computing for every $x$ a central $95\%$ interval in the normal
distribution on the right-hand side.

\begin{figure}

\includegraphics{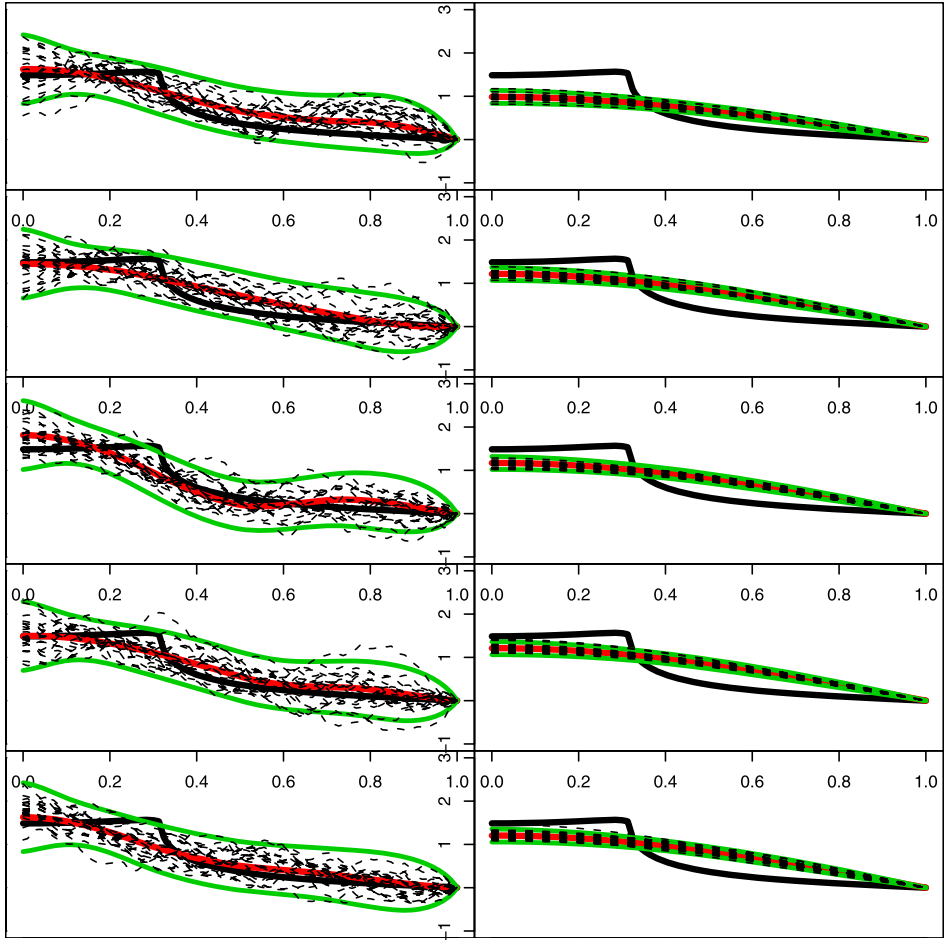}%
\vspace*{-2pt}
\caption{Realizations of the posterior mean (red) and (marginal)
posterior credible bands (green), and 20 draws from
the posterior (dashed curves). In all ten panels $n=1\mbox{,}000$ and $\b=1$.
Left 5~panels: $\a=1$; right 5 panels: $\a=5$. True curve
(black) given by coefficients $\m_{0,i}=i^{-3/2}\sin(i)$.}
\label{Figure1}
\vspace*{-2pt}
\end{figure}

Figure \ref{Figure1} illustrates these bands for
$n=1\mbox{,}000$. In every one of the 10 panels in the figure the black curve
represents the function $\m_0$, defined by the coefficients
$i^{-3/2}\sin(i)$ relative to $e_i$ ($\b=1$). The 10 panels
represent~10 independent realizations of the data, yielding 10
different realizations of the posterior mean (the red curves) and the
posterior credible bands (the green curves). In the left five
panels the prior is given by $\l_i=i^{-2\a-1}$ with $\a=1$, whereas in
the right panels the prior corresponds to $\a=5$. Each of the 10
panels also shows 20 realizations from the posterior
distribution.

%f1 ###
%fi1 #&#
%

%f2 ###
%fi2 #&#
%
\begin{figure}

\includegraphics{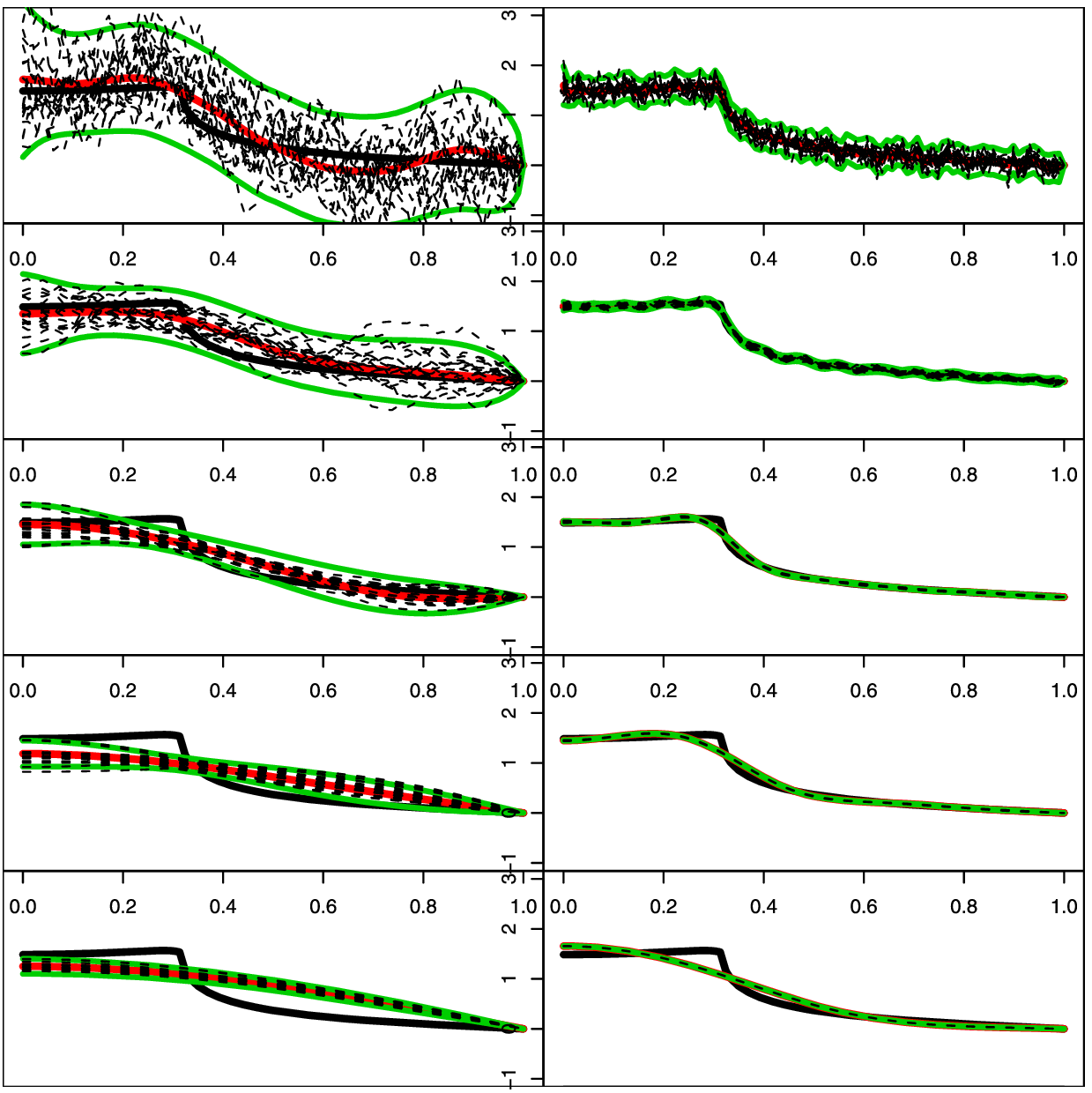}

\caption{Realizations of the posterior mean (red) and (marginal)
posterior credible bands (green), and 20 draws from
the posterior (dashed curves). In all\vspace*{1pt} ten panels $\b=1$.
Left~5 panels: $n=1\mbox{,}000$ and $\a=0.5,1,2,3,5$ (top to bottom);
right 5 panels: $n=10^8$ and $\a=0.5,1,2,3,5$ (top to bottom). True curve
(black) given by coefficients $\m_{0,i}=i^{-3/2}\sin(i)$.}
\label{Figure2}
\end{figure}

Clearly, the posterior mean is not estimating the true curve very well,
even for $n=1\mbox{,}000$. This is mostly caused by the intrinsic difficulty
of the inverse problem: better estimation requires bigger sample
size. A comparison of the left and right panels shows that
the rough prior ($\a=1$) is aware of the difficulty: it produces
credible bands that in (almost) all cases contain the true
curve. On the other hand, the smooth prior ($\a=5$)
is overconfident; the spread of the posterior distribution
poorly reflects the imprecision of estimation.

Specifying a prior that is too smooth relative to the true curve
yields a posterior distribution which gives both a bad reconstruction
and a misguided sense of uncertainty. Our theoretical results show
that the inaccurate quantification of estimation error remains even as
$n\ra\infty$.

The reconstruction, by the posterior mean or any other
posterior quantiles, will eventually converge to the true
curve. However, specification of a too smooth\vadjust{\goodbreak} prior will slow down
this convergence significantly. This is illustrated in
Figure \ref{Figure2}. Every one of its 10 panels is similarly
constructed as before, but now with $n=1\mbox{,}000$ and $n=10^8$ for the five
panels on the left-hand and right-hand side, respectively, and with
$\a=1/2,1,2,3,5$ for the five panels from top to bottom.
At first sight $\a= 1$ seems better (see the left column in Figure \ref{Figure2}),
but leads to zero coverage because of the mismatch close to the bump
(see the right column), while $\a= 1/2$ captures the bump.
For $n=10^8$ the posterior for this optimal prior has
collapsed onto the true curve, whereas the smooth posterior for $\a=5$
still has major difficulty in recovering the bump in the
true curve (even though it ``thinks'' it has captured the correct
curve, the bands having collapsed to a single curve in the\vadjust{\goodbreak}
figure).\looseness=-1

%s7 ###
%se7 #&#
\section{Proofs}\vspace*{3pt}
\label{SectionProofs}

%s7.1 ###
%su7.1 #&#
\subsection{\texorpdfstring{Proof of Theorem \protect\ref{ContrPPP}}{Proof of Theorem 4.1}}
The second moment of a Gaussian distribu\-tion on $H_1$ is equal to
the square norm of its mean plus the trace of its~co\-variance
operator.
Because the posterior is Gaussian $N(AY,S_n)$, it follows that
\[
\int\|\m-\m_0\|_1^2 \,d\Pi_n(\m\given Y)
=\|AY-\m_0\|_1^2+\tr(S_n).
\]
By Markov's inequality, the left-hand side is an upper bound to
$M_n^2{\varepsilon}_n^2\Pi_n(\m\dvtx\allowbreak \|\m-\m_0\|_1 \geq
M_n{\varepsilon}_n \given Y)$.
Therefore, it suffices to show that the expectation under~$\m_0$ of
the right-hand
side of the display is bounded by a multiple of~${\varepsilon}_n^2$.
The expectation of the
first term is the mean square error\vspace*{1pt} of the posterior mean~$AY$,
and can be written as the sum
$\|AK\m_0 - \m_0\|_1^2 +{n}^{-1}\tr(AA^T)$ of its square bias and
``variance.''
The second term $\tr(S_n)$ is deterministic. Under
Assumption \ref{common} the three quantities
can be expressed in the coefficients relative to the eigenbasis $(e_i)$
as
%
%e7.3 ###
%e7.2 ###
%e7.1 ###
%e7.1 #&#
%e7.2 #&#
%e7.3 #&#
%
\begin{eqnarray}\label{SqBias}
\|AK\m_0 - \m_0\|_1^2
&=&\sum_i\frac{\m_{0,i}^2}{(1+n\l_i\k_i^2)^2}
\asymp\sum_i \frac{\mu_{0,i}^2}{(1+n\tau_n^2i^{-1-2\alpha
-2p})^2},\\
\label{Var}
\frac{1}{n}\tr(AA^T)
&=& \sum_i\frac{ n\l_i^2\k_i^2}{(1+n\l_i\k_i^2)^2}
\asymp\sum_{i} \frac{n\tau_n^4 i^{-2-4\alpha-2p}}{(1+n\tau
_n^2i^{-1-2\alpha-2p})^2},\\
\label{PostSpr}
\tr(S_n) &=& \sum_i\frac{\l_i}{1+n\l_i\k_i^2}
\asymp\sum_{i} \frac{\tau_n^2i^{-1-2\a}}{1+n\tau_n^2i^{-1-2\alpha-2p}}.
\end{eqnarray}
By Lemma\vspace*{1pt} \ref{NormSeries} (applied with $q=\b$, $t=0$, $u = 1+2\a+2p$,
$v=2$ and $N = n\t_n^2$), the first can be bounded by\vspace*{2pt}
$\|\m_0\|^2_\b(n\t_n^2)^{-({2\b})/({1+2\a+2p})\wedge2}$, which accounts
for the first term in the definition of ${\varepsilon}_n$.
By Lemma \ref{RVSeries} [applied with
$\S(i)=1$, $q=-1/2$, $t=2+4\a+2p$, $u=1+2\a+2p$, $v=2$, and $N = n\t_n^2$],
and again Lem\-ma~\ref{RVSeries} [applied with $\S(i)=1$, $q=-1/2$,
$t=1+2\a$, $u=1+2\a+2p$, $v=1$ and $N =n\t_n^2$],
both the second and third expressions are of the order the square of the
second term in the definition of ${\varepsilon}_n$.

The consequences (i) and (ii) follow by verification after
substitution of~$\t_n$ as given. To prove consequence (iii),
we note that the two terms in the definition of ${\varepsilon}_n$ are
decreasing and increasing in $\t_n$, respectively. Therefore,
the maximum of these two terms is minimized with respect to $\t_n$ by equating
the two terms.\vspace*{1pt} This minimum (assumed at
$\t_n =n^{-({1+\a+2p})/({3+4\a+6p})}$) is much bigger than $n^{-\b
/(1+2\b+2p)}$
if $\b>1+2\a+2p$.

%s7.2 ###
%su7.2 #&#
\subsection{\texorpdfstring{Proof of Theorem \protect\ref{LinContrPPP}}{Proof of Theorem 5.1}}
By Proposition \ref{PosteriorLinear} the posterior distribution
is $N(LAY,s_n^2)$, and, hence, similarly as in the
proof of Theorem \ref{ContrPPP}, it suffices to show that
\[
\E_{\m_0} |LAY-L\m_0|^2 +s_n^2
=|LAK\m_0-L\m_0|^2+\frac1n\|LA\|_1^2+s_n^2
\]
is bounded above by a multiple of $\eps_n^2$. Under Assumption \ref{common}
the expressions on the right can be written
%
%e7.6 ###
%e7.5 ###
%e7.4 ###
%e7.4 #&#
%e7.5 #&#
%e7.6 #&#
%
\begin{eqnarray}\label{LinBias}
LAK\m_0-L\m_0
&=& -\sum_{i} \frac{l_i\m_{0,i}}{1+n\l_i\k_i^2}
\lesssim\sum_i \frac{|l_i\m_{0,i}|}{1+n\t_n^2i^{-1-2\a-2p}},\\
\label{LinVar}
t_n^2:\!&=&\frac{1}{n}\|LA\|_1^2
= \sum_{i}
\frac{l_i^2n\l_i^2\k_i^2}{(1+n\l_i\k_i^2)^2}\nonumber\\[-9pt]\\[-9pt]
&\asymp& n\t_n^4 \sum_i \frac{l_i^2i^{-2-4\a-2p}}{(1+n\t
_n^2i^{-1-2\a-2p})^2},\nonumber\\[-2pt]
\label{LinPostSpr}
s_n^2 &=& \sum_{i} \frac{l_i^2\l_i}{1+n\l_i\k_i^2}
\asymp\t_n^2 \sum_{i} \frac{l_i^2 i^{-1-2\a}}{1+n\t_n^2i^{-1-2\a-2p}}.
\end{eqnarray}
By the Cauchy--Schwarz inequality the square of the bias
(\ref{LinBias}) satisfies
%
%e7.7 ###
%e7.7 #&#
%
\begin{equation}
\label{EqSqBiasLinear}
|LAK\m_0-L\m_0|^2
\lesssim\|\m_0\|^2_\b\sum_i \frac{l_i^2i^{-2\b}}{(1+n\t
_n^2i^{-1-2\a-2p})^2}.
\end{equation}
By Lemma \ref{NormSeries}
(applied with $q=q, t=2\b, u=1+2\a+2p, v = 2$ and $N = n\t_n^2$)
the right-hand side of this display
can be further bounded by $\|\m_0\|_\b^2\|l\|_q^2$ times the square of
the first term in the sum of two terms that defines ${\varepsilon}_n$.
By Lemma \ref{NormSeries} (applied
with $q=q, t=2+4\a+2p, u = 1+2\a+2p, v=2$ and $ N = n\t_n^2$)
and again Lemma \ref{NormSeries} (applied with $q=q, t=1+2\a, u=1+2\a
+2p, v = 1$ and $N = n\t_n^2$), the right-hand sides of (\ref{LinVar})
and (\ref{LinPostSpr})
are bounded above by $\|l\|^2_q$ times the square of the second term
in the definition of ${\varepsilon}_n$.\looseness=-1

Consequences (i)--(iv) follow by substitution, and, in the
case of (iii), optimization over $\t_n$.\vspace*{-3pt}

%s7.3 ###
%su7.3 #&#
\subsection{\texorpdfstring{Proof of Theorem \protect\ref{LinContrPPPRV}}{Proof of Theorem 5.2}}
This follows the same lines as the proof of Theorem \ref{LinContrPPP},
except that we use Lemma \ref{RVSeries}
(with $q\,{=}\,q, t\,{=}\,2\b$, \mbox{$u\,{=}\,1\,{+}\,2\a\,{+}\,2p$}, \mbox{$v\,{=}\,2$} and $N\,{=}\,n\t_n^2$)
and Lemma \ref{RVSeries}
(with $q\,{=}\,q, t\,{=}\,2\,{+}\,4\a\,{+}\,2p, u\,{=}\,1\,{+}\,2\a\,{+}\,2p$, \mbox{$v=2$} and $N = n\t_n^2$)
and again Lemma \ref{RVSeries}
(with $q=q, t=1+2\a, u=1+2\a+2p$, \mbox{$v = 1$} and $N = n\t_n^2$)
to bound the three terms
(\ref{LinVar})--(\ref{EqSqBiasLinear}).\vspace*{-3pt}

%s7.4 ###
%su7.4 #&#
\subsection{\texorpdfstring{Proof of Theorem \protect\ref{CrNS}}{Proof of Theorem 4.2}}
Because the posterior distribution is $N(AY,\allowbreak S_n)$,
by Proposition \ref{Posterior}, the radius\vspace*{1pt} $r_{n,\g}$ in
(\ref{EqRadius}) satisfies
$\Pr(U_n< r_{n,\g}^2)=1-\g$,
for~$U_n$ a random variable distributed as the square norm
of an $N(0,S_n)$-variable. Under (\ref{EqProblem}) the variable
$AY$ is $N(AK\m_0,n^{-1}AA^T)$-distributed, and, thus, the coverage
(\ref{EqCoverage}) can be written as
%
%e7.8 ###
%e7.8 #&#
%
\begin{equation}
\label{EqCoverageW}
\Pr(\|W_n+AK\m_0-\m_0\|_1\le r_{n,\g})
\end{equation}
for $W_n$ possessing\vspace*{1pt} a $N(0,n^{-1}AA^T)$-distribution.
For ease of notation let $V_n=\|W_n\|_1^2$.

The variables $U_n$ and $V_n$ can be represented as
$U_n=\sum_i s_{i,n}Z_i^2$ and $V_n=\sum_i t_{i,n}Z_i^2$,
for $Z_1,Z_2,\ldots$ independent standard normal variables,
and $s_{i,n}$ and $t_{i,n}$ the eigenvalues of $S_n$ and $n^{-1}AA^T$,
respectively, which satisfy
\begin{eqnarray*}
s_{i,n}&=&\frac{\l_i}{1+n\l_i\k_i^2}
\asymp\frac{\t_n^2i^{-2\a-1}}{1+n\t_n^2i^{-2\a-2p-1}},\\
t_{i,n}&=&\frac{ n\l_i^2\k_i^2}{(1+n\l_i\k_i^2)^2}
\asymp\frac{n\t_n^4i^{-4\a-2p-2}}{(1+n\t_n^2i^{-2\a-2p-1})^2},\\
s_{i,n}-t_{i,n}&=&\frac{ \l_i}{(1+n\l_i\k_i^2)^2}
\asymp\frac{\t_n^2i^{-2\a-1}}{(1+n\t_n^2i^{-2\a-2p-1})^2}.
\end{eqnarray*}
Therefore, by Lemma \ref{RVSeries} (applied with $\S\equiv1$ and $q=-1/2$;
always the first case),
\begin{eqnarray*}
\E U_n&=&\sum_i s_{i,n}\asymp\t_n^2(n\t_n^2)^{-{2\a}/({1+2\a
+2p})},\\
\E V_n&=&\sum_i t_{i,n}\asymp\t_n^2(n\t_n^2)^{-{2\a}/({1+2\a
+2p})},\\
\E(U_n-V_n)&=&\sum_i(s_{i,n}-t_{i,n})\asymp\t_n^2(n\t_n^2)^{-
{2\a}/({1+2\a+2p})},\\
\var U_n&=&2\sum_i s_{i,n}^2\asymp\t_n^4(n\t_n^2)^{-({1+4\a
})/({1+2\a+2p})},\\
\var V_n&=&2\sum_i t_{i,n}^2\asymp\t_n^4(n\t_n^2)^{-({1+4\a
})/({1+2\a+2p})}.
\end{eqnarray*}
We conclude that the standard deviations of $U_n$ and $V_n$ are
negligible relative to their means, and also relative to the
difference $\E(U_n-V_n)$ of their means. Because $U_n\ge V_n$, we
conclude that the distributions of $U_n$ and $V_n$ are asymptotically
completely separated: $\Pr(V_n\le v_n\le U_n)\to1$ for some $v_n$
[e.g., $v_n=\E(U_n+V_n)/2$]. The numbers $r_{n,\g}^2$ are
$1-\g$-quantiles of $U_n$, and, hence, \mbox{$\Pr(V_n\le r^2_{n,\g
}(1+o(1)))\ra
1$}. Furthermore, it follows that
\[
r_{n,\g}^2 \asymp\t_n^2(n\t_n^2)^{-{2\a}/({1+2\a+2p})}\asymp
\E U_n\asymp\E V_n.
\]
The square norm of the bias $AK\m_0-\m_0$ is given in
(\ref{SqBias}), where it was noted that
\[
B_n:={\sup_{\|\m_0\|_\b\lesssim1}}\|AK\m_0-\m_0\|_1 \asymp
(n\t_n^2)^{-\b/({1+2\a+2p})\wedge1}.
\]
The bias $B_n$ is decreasing
in $\t_n$, whereas $\E U_n$ and $\var U_n$ are increasing. The scaling
rate $\tilde\t_n\asymp n^{(\a- \wt\b)/(1+2\wt\b+2p)}$ balances the
square bias $B_n^2$ with the variance $\E V_n$ of the posterior mean,
and hence with $r_{n,\g}^2$.

Case (i). In this case $B_n \ll r_{n,\g}$. Hence,
$\Pr(\|W_n+AK\m_0-\m_0\|_1\le r_{n,\g})
\ge\Pr(\|W_n\|_1\le r_{n,\g}-B_n)
=\Pr(V_n\le r_{n,\g}^2(1+o(1)))\ra1$,
uniformly in the set of~$\m_0$ in the supremum defining $B_n$.

Case (iii). In this case $B_n \gg r_{n,\g}$. Hence,
$\Pr(\|W_n+AK\m_0^n-\m_0^n\|_1\le r_{n,\g}) \le
\Pr(\|W_n\|_1\ge B_n-r_{n,\g})\ra0$ for any sequence
$\m_0^n$ (nearly) attaining the supremum in the definition of
$B_n$. If $\t_n\equiv1$, then $B_n$ and $r_{n,\g}$ are both
powers of $1/n$ and, hence, $B_n \gg r_{n,\g}$ implies that
$B_n\gtrsim r_{n,\g} n^\d$, for some $\d>0$. The preceding argument
then applies for a fixed $\m_0$ of the form $\m_{0,i}\asymp
i^{-1/2-\b-{\varepsilon}}$, for small ${\varepsilon}>0$, that gives
a bias that is much
closer than $n^\d$ to $B_n$.

Case (ii). In this case $B_n\asymp r_{n,\g}$. If $\b<1+2\a+2p$, then by
the second assertion (first case) of Lemma \ref{NormSeries} the bias
$\|AK\m_0-\m_0\|_1$ at a fixed $\m_0$ is of strictly smaller order
than the supremum $B_n$. The argument of (i) shows that the asymptotic
coverage then tends to 1.

Finally, we prove the existence of a sequence $\m_0^n$ along which
the coverage is a given $c\in[0,1)$. The coverage (\ref{EqCoverageW})
with $\m_0$ replaced by $\m_0^n$ tends to $c$ if,
for $b_n=AK\m_0^n-\m_0^n$ and $z_c$ a standard normal quantile,
%
%e7.10 ###
%e7.9 ###
%e7.9 #&#
%e7.10 #&#
%
\begin{eqnarray}
\label{EqNormality}
\frac{\|W_n+b_n\|_1^2-\E\|W_n+b_n\|_1^2}{{\sdev}\|W_n+b_n\|_1^2}
&\weak& N(0,1),\\
\label{EqConvergenceRng}
\frac{r_{n,\g}^2-\E\|W_n+b_n\|_1^2}{{\sdev}\|W_n+b_n\|_1^2}
&\ra& z_c.
\end{eqnarray}
Because $W_n$ is mean-zero Gaussian, we have
$\E\|W_n+b_n\|_1^2=\E\|W_n\|_1^2+\|b_n\|_1^2$ and
${\var}\|W_n+b_n\|_1^2={\var}\|W_n\|_1^2+4\var\langle W_n,b_n\rangle_1$.
Here $\|W_n\|_1^2=V_n$ and the distribution of
$\langle W_n,b_n\rangle_1$ is zero-mean Gaussian with variance
$\langle b_n, n^{-1}AA^Tb_n\rangle_1$. With
$t_{i,n}$ the eigenvalues of $n^{-1}AA^T$, display (\ref
{EqConvergenceRng}) can be
translated in the coefficients $(b_{n,i})$ of $b_n$ relative
to the eigenbasis, as
%
%e7.11 ###
%e7.11 #&#
%
\begin{equation}
\label{Hulp}
\frac{r_{n,\g}^2-\E V_n-\sum_i b_{n,i}^2}{\sqrt{\var V_n+4\sum
_it_{i,n}b_{n,i}^2}}
\ra z_c.
\end{equation}
We choose $(b_{n,i})$ differently in the cases that $\b\le1+2\a+2p$
and $\b\ge1+2\a+2p$, respectively. In both cases the sequence
has exactly one nonzero coordinate. We denote this coordinate by $b_{n,i_n}$,
and set, for numbers $d_n$ to be determined,
\[
b_{n,i_n}^2=r_{n,\g}^2-\E V_n-d_n\sdev V_n.
\]
Because $r_{n,\g}^2$, $\E V_n$ and $r_{n,\g}^2-\E V_n$ are of the
same order
of magnitude, and~$\sdev V_n$ is of strictly smaller
order, for bounded or slowly diverging $d_n$ the right-hand side of the
preceding display is equivalent to $(r_{n,\g}^2-\E V_n)(1+o(1))$.
Consequently, the left-hand side of (\ref{Hulp}) is equivalent to
\[
\frac{d_n \sdev V_n}{\sqrt{\var V_n+4t_{i_n,n}(r_{n,\g}^2-\E V_n)(1+o(1))}}.
\]
The remainder of the argument is different in the two cases.

Case $\b\le1+2\a+2p$. We choose $i_n\asymp(n\t_n^2)^{1/(1+2\a+2p)}$.
It can be verified that $t_{i_n,n}(r_{n,\g}^2-\E V_n)/\var V_n\asymp1$.
Therefore, for $c\in[0,1]$, there exists a bounded or slowly
diverging sequence $d_n$
such that the preceding display tends to $z_c$.\vadjust{\goodbreak}

The bias $b_n$ results from a parameter $\m_0^n$ such that
$b_{n,i}=(1+n \l_i\k_i^2)^{-1}(\m_0^n)_i$, for every $i$. Thus,
$\m_0^n$ also has exactly one nonzero coordinate, and this is
proportional to the corresponding coordinate of $b_n$, by the definition
of $i_n$. It follows that
\[
i_n^{2\b}(\m_0^n)^2_{i_n}\asymp i_n^{2\b} b_{n,i_n}^2
\lesssim i_n^{2\b}(r_{n,\g}^2-\E V_n)\asymp1
\]
by the definition of $\t_n$. It follows that $\|\m_0^n\|_\b\lesssim1$.

Case $\b\ge1+2\a+2p$. We choose $i_n=1$. In this case
$\t_n\ra0$ and it can be verified that
$t_{i_n,n}(r_{n,\g}^2-\E V_n)/\var V_n\ra0$. Also,
\[
(\m_0^n)_1^2\asymp(1+n\t_n^2)^2b_{n,1}^2\lesssim(1+n\t_n^2)^2\E V_n.
\]
This is $O(1)$, because $\t_n$ is chosen so that $\E V_n$
is of the same order as the square bias $B_n^2$, which is $(n\t_n^2)^{-2}$
in this case.

It remains to prove the asymptotic normality (\ref{EqNormality}).
We can write
\[
\|W_n+b_n\|_1^2-\E\|W_n+b_n\|_1^2
= \sum_i t_{i,n}(Z_i^2-1)+2b_{n,i_n}\sqrt{t_{{i_n},n}}Z_{i_n}.
\]
The second term is normal by construction. The first
term has variance $2\sum_i t_{i,n}^2$. With some effort it can be seen that
\[
\sup_i \frac{t_{i,n}^2}{\sum_i t_{i,n}^2}\ra0.
\]
Therefore, by a slight adaptation of the Lindeberg--Feller theorem
(to infinite sums), we have that $\sum_i t_{i,n}(Z_i^2-1)$ divided
by its standard deviation tends in distribution to the standard
normal distribution. Furthermore, the preceding display shows
that this conclusion does not change if the $i_n$th term is left out
from the infinite sum. Thus, the two terms converge jointly
to asymptotically independent standard normal variables, if scaled
separately by their standard deviations.
Then their scaled sum is also asymptotically
standard normally distributed.

%s7.5 ###
%su7.5 #&#
\subsection{\texorpdfstring{Proof of Theorem \protect\ref{LinCrS}}{Proof of Theorem 5.3}}
Under (\ref{EqProblem}) the variable $LAY$ is
$N(LAK\m_0,\allowbreak t_n^2)$-distributed, for $t_n^2$ given in (\ref{LinVar}).
It follows that the coverage can be written, with $W$
a standard normal variable,
%
%e7.12 ###
%e7.12 #&#
%
\begin{equation}
\label{EqCoverageLinW}
\Pr( |Wt_n+LAK\m_0-L\m_0|\le-s_n z_{\g/2}).
\end{equation}
The bias $LAK\m_0-L\m_0$ and posterior spread $s_n^2$ are
expressed as a series in~(\ref{LinBias}) and (\ref{LinPostSpr}).

In the proof of Theorem \ref{LinContrPPPRV}
$s_n$ and $t_n$ were seen to have the same order of magnitude, given by
the second term in ${\varepsilon}_n$ given in (\ref{EqTauLRV}),
with a slowly varying term $\d_n$ as given in the theorem,
%
%e7.13 ###
%e7.13 #&#
%
\begin{equation}
\label{EqSnTnLin}
s_n\asymp t_n\asymp\t_n (n\t_n^2)^{-(1/2+\a+q)/(1+2\a+2p)}\d_n.
\end{equation}
Furthermore,\vspace*{1pt} $t_n \leq s_n$ for every $n$, as every term in the
infinite series (\ref{LinVar}) is $n\l_i\k_i^2/(1+n\l_i\k_i^2)\le1$
times the corresponding term in (\ref{LinPostSpr}).

Because $W$ is centered, the coverage (\ref{EqCoverageLinW}) is
largest if the bias $LAK\m_0-L\m_0$ is zero. It is then at least
$1-\g$, because $t_n\le s_n$; remains strictly smaller than $1$,
because $t_n
\asymp s_n$; and tends to exactly $1-\g$ iff $s_n/t_n \ra1$. By
Theorem~\ref{TBvM-reg}(i) the latter is impossible if $q < p$.
The analysis for nonzero $\m_0$ depends strongly on the size
of the bias relative to $t_n$.

The supremum of the bias satisfies, for $\g_n$ the slowly varying term
given in Theorem \ref{LinContrPPPRV},
%
%e7.14 ###
%e7.14 #&#
%
\begin{equation}
\label{EqSupBiasLinear}
B_n:={\sup_{\|\m_0\|_\b\lesssim1}}|LAK\m_0-L\m_0|
\asymp(n\t_n^2)^{-((\b+q)/(1+2\a+2p))\wedge1}\g_n.
\end{equation}
That the left-hand side of (\ref{EqSupBiasLinear})
is smaller than the right-hand side was already shown
in the proof of Theorem \ref{LinContrPPPRV}, with
the help of Lemma \ref{RVSeries}. That this upper bound is
sharp follows by considering the sequence $\m_0^n$ defined by,
with $\tilde B_n$ the right-hand side of the preceding display,
\[
\m_{0,i}^{n} = \frac1{\tilde B_n} \frac{i^{-2\b}l_i}{1+n\t
_n^2i^{-1-2\a-2p}}.
\]
[This is the sequence that gives equality in the application
of the Cauchy--Schwarz inequality to derive (\ref{EqSqBiasLinear}).]
Using Lemma \ref{RVSeries}, it can be seen that $\|\m_0^n\|_\b
\lesssim1$
and that the bias at $\m_0^n$ is of the order $\tilde B_n$.

By Lemma \ref{LemmaTechnicalBias},
the bias at a \textit{fixed} $\m_0\in S^\b$ is of strictly smaller order
than the supremum $B_n$ if $\b+q<1+2\a+2p$.

The maximal bias $B_n$ is a decreasing function of the scaling
parameter~$\t_n$, while the standard deviation $t_n$ and root-spread
$s_n$ increase with $\t_n$.
The scaling rate $\wt\t_n$ in the statement of the
theorem balances $B_n$ with $s_n\asymp t_n$.

Case (i). If $\t_n\gg\wt\t_n$, then $B_n\ll t_n$. Hence, the bias
$LAK\m_0-L\m_0$ in (\ref{EqCoverageLinW}) is negligible relative
to $t_n\asymp s_n$, uniformly in $\|\m_0\|_\b\lesssim1$, and
the coverage is asymptotic to $\Pr( |Wt_n|\le-s_nz_{\g/2})$,
which is asymptotically strictly between $1-\g$ and~$1$.

Case (iii). If $\t_n\ll\wt\t_n$, then $B_n\gg t_n$.
If $b_n=LAK\m_0^n-L\m_0^n$ is
the bias at a~sequence $\m_0^n$ that (nearly) attains the supremum
in the definition of~$B_n$, then the coverage at $\m_0^n$ satisfies
$\Pr( |Wt_n+b_n|\le-s_n z_{\g/2})
\le\Pr( |Wt_n|\ge b_n-s_n |z_{\g/2}|)\ra0$,
as $b_n\asymp B_n\gg s_n$. By the same argument,
the coverage also tends to zero for a fixed $\m_0$
in $S^\b$ with bias $b_n=LAK\m_0-L\m_0\gg t_n$.
For this we choose $\m_{0,i}=i^{-\b}i^ql_i\tilde S(i)$ for a slowly
varying function such that $\sum_i \S^2(i)\tilde S^2(i)/i<\infty$.
The latter condition ensures that $\|\m_0\|_\b<\infty$.
By another application of Lemma \ref{RVSeries}, the bias
at $\m_0$ is of the order [cf. (\ref{LinBias})]
\[
\sum_i\frac{l_i\m_{0,i}}{1+n\t_n^2i^{-1-2\a-2p}}
=\sum_i \frac{(l_i\tilde S^{1/2}(i))^2i^{-\b+q}}{1+n\t_n^2i^{-1-2\a-2p}}
\asymp(n\t_n^2)^{-({\b+q})/({1+2\a+2p})\wedge1}\tilde\g_n,
\]
where, for $\r_n=(n\t_n^2)^{1/(1+2\a+2p)}$,
\[
\tilde\g_n^2=\cases{
\S^2(\r_n)\tilde S(\r_n), &\quad if $\b+q < 1+2\a+2p$,\vspace*{2pt}\cr
\displaystyle \sum_{i \leq\r_n}\frac{\S^2(i)\tilde S(i)}{i}, &\quad if
$\b+q=1+2\a+2p$,\vspace*{2pt}\cr
1, &\quad if $\b+q>1+2\a+2p$.}
\]
Therefore, the bias at $\m_0$ has the same form as the maximal
bias $B_n$; the difference is in the slowly varying factor $\tilde\g_n$.
If $\t_n\le\tilde\t_nn^{-\d}$, then $B_n\gtrsim t_nn^{\d'}$ for
some $\d'>0$
and, hence, $b_n\asymp B_n\tilde\g_n/\g_n\gg t_n$.

Case (ii). If $\t_n\asymp\wt\t_n$, then $B_n\asymp t_n$. If
$b_n=LAK\m_0^n-L\m_0^n$ is again the bias at a sequence $\m_0^n$ that
nearly assumes the supremum in the definition of $B_n$, we have that
$\Pr( |Wt_n+d b_n|\le-s_n z_{\g/2}) \le\Pr(
|Wt_n|\ge
db_n-s_n |z_{\g/2}|)$ attains an arbitrarily small value if $d$
is chosen sufficiently large. This is the coverage
at the sequence~$d\m_0^n$, which is bounded in $S^\b$.
On the other hand, the bias at a fixed
$\m_0\in S^\b$ is of strictly smaller order than the supremum~$B_n$,
and, hence, the coverage at a fixed $\m_0$ is as in case (i).

If the scaling rate is fixed to $\t_n\equiv1$, then it can be checked
from (\ref{EqSnTnLin}) and~(\ref{EqSupBiasLinear}) that $B_n\ll t_n$,
$B_n\asymp t_n$ and $B_n\gg t_n$ in the three cases $\a<\b-1/2$,
$\a=\b-1/2$ and $\a>\b-1/2$, respectively. In the first and third
cases the maximal bias and the spread differ by more than a polynomial term
$n^\d$; in the second case it must be noted that the slowly varying
terms $\g_n$ and $\d_n$ are equal [to $\S(\r_n)$]. It follows that the
preceding analysis (i), (ii), (iii) extends to this situation.

%s7.6 ###
%su7.6 #&#
\subsection{\texorpdfstring{Proof of Theorem \protect\ref{TBvM-reg}}{Proof of Theorem 5.4}}
(i). The two quantities $s_n$ and $t_n$ are given as series in
(\ref{LinPostSpr}) and (\ref{LinVar}). Every term in the series
(\ref{LinVar}) is $n\l_i\k_i^2/(1+n\l_i\k_i^2)\le1$ times the
corresponding term in the series (\ref{LinPostSpr}). Therefore,
$s_n/t_n\ra1$ if and only if the series are determined by the terms for
which these numbers are ``close to'' 1, that is, $n\l_i\k_i^2$ is
large. More precisely, we show below that $s_n/t_n\ra1$ if and only
if, for every $c > 0$,
%
%e7.15 ###
%e7.15 #&#
%
\begin{equation}
\label{Hulp1}
\sum_{n\l_i\k_i^2 \leq c} \frac{l_i^2\l_i}{1+n\l_i\k_i^2}
= o\biggl(\sum_i \frac{l_i^2\l_i}{1+n\l_i\k_i^2}\biggr).
\end{equation}
If $l\in S^p$, then the series on the left is as in
Lemma \ref{NormSeries}
with $q=p$, $u=1+2\a+2p$, $v=1$, $N=n\t_n^2$ and $t=1+2\a$. Hence,
$(t+2q)/u\ge v$, and the display follows from the final assertion
of the lemma. If $l_i=i^{-q-1/2}\S(i)$ for a slowly varying function
$\S$,
then the series is as in Lemma \ref{RVSeries}, with the
same parameters, and by the last statement
of the lemma the display is true if and only if
$(t+2q)/u\ge v$, that is, $q\ge p$.

To prove that (\ref{Hulp1}) holds iff $s_n/t_n\ra1$,
write $s_n^2=A_n+B_n$, for $A_n$ and $B_n$ the
sums over the terms in (\ref{LinPostSpr})
with $n\l_i\k_i^2> c$ and $n\l_i\k_i^2\le c$,
respectively, and, similarly, $t_n^2=C_n+D_n$. Then
\[
\frac{D_n}{B_n}\le\frac c{1+c}\le\frac{C_n}{A_n}\le1.
\]
It follows that
\[
\frac{t_n^2}{s_n^2}
=\frac{C_n+D_n}{A_n+B_n}
=\frac{C_n/A_n+(D_n/B_n)(B_n/A_n)}{1+B_n/A_n}
\le\frac{1+c/(1+c)(B_n/A_n)}{1+B_n/A_n}.
\]
Because $x\mapsto(1+rx)/(1+x)$ is strictly decreasing from 1 at $x=0$
to $r<1$ at $x=\infty$ (if $0<r<1$), the right-hand side
of the equation is asymptotically~1 if and only if $B_n/A_n\to0$,
and otherwise its liminf is strictly smaller.
Thus, $t_n/s_n\to1$ implies that $B_n/A_n\to0$.
Second,
\[
\frac{t_n^2}{s_n^2}
\ge\frac{C_n}{A_n+B_n}
=\frac{C_n/A_n}{1+B_n/A_n}\ge\frac{c/(1+c)}{1+B_n/A_n}.
\]
It follows that $\liminf t_n^2/s_n^2\ge c/(1+c)$ if $B_n/A_n\to0$.
This being true for every $c>0$ implies that $t_n/s_n\to1$.

\hphantom{i}(i) Second assertion.\vspace*{1pt} If $l\in S^p$, then we apply
Lemma \ref{NormSeries} with $q=p$, $t=1+2\a$, $u=1+2\a+2p$, $v=1$ and
$N=n\t_n^2$ to see that $s_n^2\asymp\t_n^2(n\t_n^2)^{-v}=n^{-1}$.
Furthermore, the second assertion of the lemma with $(uv-t)/2=p$ shows
that $ns_n^2\ra\|l\|_{p}^2=\sum_il_i^2/\k_i^2$ in the case that
$\k_i=i^{-p}$. The proof can be extended to cover the slightly more
general sequence $(\k_i)$ in Assumption~\ref{PPP}.

If $l\in\RV^q$, then we apply Lemma \ref{RVSeries} with $q=p$,
$t=1+2\a$, $u=1+2\a+2p$, $v=1$ and $N=n\t_n^2$ to see that
$s_n^2\asymp n^{-1}\sum_{i\le N^{1/u}}\S^2(i)/i$.\vspace*{1pt}

(ii) If $l\in S^q$, then the bias is bounded above in
(\ref{EqSqBiasLinear}), and in the proof of Theorem \ref{LinContrPPP} its
supremum $B_n$ over $\|\m_0\|_\b\lesssim1$ is seen to be bounded by
$(n\t_n^2)^{-(\b+q)/(1+2\a+2p)\wedge1}$, the first term in the
definition of ${\varepsilon}_n$ in the statement of this theorem. This
upper bound
is $o(n^{-1/2})$ iff the stated conditions hold.
[Here we use that $\S^2(N)\ll\sum_{i\le N}\S^2(i)/i$
as $N\ra\infty$, as noted in the proof of Lem\-ma~\ref{RVSeries}.]

The supremum of the bias $B_n$ in the case that $l\in\RV^q$ is given in
(\ref{EqSupBiasLinear}). It was already seen to be $o(t_n)$ if
$\t\gg\wt\t_n$ in the proof of case\vspace*{1pt} (i) of Theorem~\ref{LinCrS}.
If $\t_n=1$, we have that $B_n\asymp n^{-(\b+q)/(1+2\a+2p)\wedge
1}\g_n$,
for $\g_n$ the slowly varying factor
given in the statement of Theorem \ref{LinContrPPPRV}. Furthermore,
we have $s_n\asymp t_n\asymp n^{-1/2}\d_n$, for $\d_n$ the slowly varying
factor in the same statement. Under the present conditions,
\mbox{$\d_n\asymp1$} if $q>p$ and $\d_n^2\asymp\sum_{i\le\r_n}\S^2(i)/i$
if $q=p$. We can now verify that $B_n=o(t_n)$ if and only if the
conditions as stated hold.

(iii) The total variation distance between two Gaussian distributions
with the same expectation and standard deviations $s_n$ and $t_n$
tends to zero if and only if $s_n/t_n\ra1$. Similarly, the total
distance between two Gaussians with the same standard deviation $s_n$
and means $\m_n$ and $\n_n$ tends to zero if and only if
$\m_n-\n_n=o(s_n)$. Therefore, it suffices to show that $(LAY-\sum_i
Y_il_i/\k_i)/s_n\ra0$ if $l\in S^p$. Because the bias was already
seen to be $o(t_n)$ and $s_n\asymp n^{-1/2}$ if $l\in S^p$, it
suffices to show that $LAZ-\sum_i Z_il_i/\k_i\ra0$.
Under Assumption \ref{PPP} this difference is equal to
\[
\sum_i \frac{\k_i\l_il_iZ_i}{n^{-1}+\k_i^2\l_i}-\sum_iZ_i\frac
{l_i}{\k_i}
=\sum_i \frac{Z_il_i}{\k_i}\biggl(\frac1{1+n\k_i^2\l_i}\biggr).
\]
If $\sum_il_i^2/\k_i^2<\infty$, then the variance of this
expression is seen to tend to zero by dominated convergence.

The final assertion of the theorem follows along the
lines of the proof of Theorem \ref{LinCrS}.

%s8 ###
%se8 #&#
\section{Technical lemmas}
\label{SectionTechnical}

%le8.1 #&#
%
\begin{lemma}\label{NormSeries}
For any $q \geq0$, $t \geq-2q$, $u >0$ and $ v\ge0$,
as $N \to\infty$,
\[
\sup_{\|\xi\|_q\le1}\sum_i \frac{\xi_i^2 i^{-t}}{(1+N i^{-u})^v}
\asymp N^{-((t+2q)/u)\wedge v}.
\]
Moreover, for every fixed $\xi\in S^q$, as $N \to\infty$,
\[
N^{((t+2q)/u)\wedge v}\sum_i \frac{\xi_i^2 i^{-t}}{(1+N i^{-u})^v}
\ra
\cases{
0, &\quad if $(t+2q)/u < v$,\cr
\|\xi\|^2_{(uv-t)/{2}}, &\quad if $(t+2q)/u \ge v$.}
\]
The last assertion remains true if the sum is limited
to the terms $i\le c N^{1/u}$, for any $c>0$.
\end{lemma}
\begin{pf}
In the range $i\le N^{1/u}$ we have $Ni^{-u}\le1+N i^{-u}\le2Ni^{-u}$,
while $1\le1+N i^{-u}\le2$ in the range $i>N^{1/u}$.
Thus, deleting either the first or second term, we obtain
\begin{eqnarray*}
\sum_{i \leq N^{1/u}} \frac{\xi_i^2 i^{-t}}{(1+N i^{-u})^v}
&\asymp&\sum_{i \leq N^{1/u}} \xi_i^2 i^{2q} \frac{i^{uv-t-2q}}{N^v}
\leq\|\xi\|_q^2 N^{-((t+2q)/u)\wedge v},\\
\sum_{i > N^{1/u}} \frac{\xi_i^2 i^{-t}}{(1+N i^{-u})^v}
&\asymp&\sum_{i > N^{1/u}} \xi_i^2i^{2q} i^{-t-2q}
\leq N^{-(t+2q)/u} \sum_{i > N^{1/u}}\xi_i^2i^{2q}.
\end{eqnarray*}
The inequality in the first line follows by bounding $i$ in $i^{uv-t-2q}$
by $N^{1/u}$ if $uv-t-2q>0$, and by 1 otherwise.
This proves the upper bound for the supremum.

The lower bound follows by considering the two sequences $(\xi_i)$
given by $\xi_i=i^{-q}$ for $i\sim N^{1/u}$ and $\xi_i=0$ otherwise
(showing that the supremum is bigger than $N^{-(t+2q)/u}$), and given
by $\xi_1=1$ and $\xi_i=0$ otherwise (showing that the supremum is
bigger than $N^{-v}$).

The second line of the preceding display shows that the sum over the
terms $i>N^{1/u}$ is $o(N^{-(t+2q)/u})$. Furthermore, the first line
can be multiplied by $N^{(t+2q)/u}$ to obtain
\[
N^{(t+2q)/u}\sum_{i \leq N^{1/u}} \frac{\xi_i^2 i^{-t}}{(1+N i^{-u})^v}
\asymp\sum_{i \leq N^{1/u}}\xi_i^2 i^{2q} \biggl(\frac
{i}{N^{1/u}}\biggr)^{uv-t-2q}.
\]
If $(t+2q)/u<v$, then $uv-t-2q>0$ and this tends to
zero by dominated convergence. Also,
\[
N^v\sum_i \frac{\xi_i^2 i^{-t}}{(1+N i^{-u})^v}
=\sum_i \xi_i^2i^{uv-t} \biggl(\frac{ N i^{-u}}{1+N i^{-u}}\biggr)^v.
\]
If $(t+2q)/u\ge v$, then $q\ge(uv-t)/2$ and, hence, $\xi\in
S^{{(uv-t)}/{2}}$, and
the right-hand side tends to $\sum_i \xi_i^2i^{uv-t}$ by dominated
convergence.

The final assertion needs to be proved only in the case that
$(t+2q)/u\ge v$,
as in the other case the whole sum tends to 0.
The sum over the terms $i>N^{1/u}$ was seen to be always
$o(N^{-(t+2q)/u})$, which is $o(N^{-v})$ if $(t+2q)/u\ge v$. The final
assertion for $c=1$ follows, because the sum over the terms $i\le N^{1/u}$
was seen to have the exact order $N^{-v}$ (if $\xi\not=0$). For general
$c$ the proof is analogous, or follows by scaling $N$.
\end{pf}
%
%le8.2 #&#
%
\begin{lemma}\label{RVSeries}
For any $t,v\ge0$, $u>0$, and $(\xi_i)$ such that
$|\xi_i|=i^{-q-1/2}\S(i)$ for $q > -t/2$ and a slowly varying function
$\S\dvtx (0,\infty)\to(0,\infty)$, as $N \ra\infty$,
\[
\sum_i \frac{\xi_i^2i^{-t}}{(1+N i^{-u})^v}
\asymp\cases{ N^{-(t+2q)/u}\S^2(N^{1/u}), &\quad if $(t+2q)/u<v$,\vspace*{2pt}\cr
\displaystyle N^{-v}\sum_{i \leq N^{1/u}}\S^2(i)\big/i, &\quad if $(t+2q)/u=v$,\vspace*{2pt}\cr
N^{-v}, &\quad if $(t+2q)/u>v$.}
\]
Moreover, for every $c>0$, the sum on the left
is asymptotically equivalent to the same sum restricted to the
terms $i\le cN^{1/u}$ if and only if $(t+2q)/u\ge v$.
\end{lemma}
\begin{pf}
As in the proof of the preceding lemma, we split the
infinite series in the sum over the terms $i\le N^{1/u}$ and $i>N^{1/u}$.
For the first part of the series
\[
\sum_{i \leq N^{1/u}} \frac{\xi_i^2 i^{-t}}{(1+N i^{-u})^v}
\asymp\sum_{i \leq N^{1/u}} \S(i)^2 \frac{i^{uv-t-2q-1}}{N^v}.
\]
If $uv-t-2q>0$ [i.e., $(t+2q)/u<v$], the right-hand side
is of the order $N^{-(t+2q)/u}\S^2(N^{1/u})$,
by Theorem 1(b) on page 281 in \cite{Feller},
while if $uv-t-2q<0$, it is of the order $N^{-v}$
by Lemma on page 280 in \cite{Feller}.
Finally, if $uv-t-2q=0$, then the right-hand side is identical to
$N^{-v}\sum_{i \leq N^{1/u}} {\S^2(i)}/{i}$.

The other part of the infinite series satisfies,
by Theorem 1(a) on page~281 in \cite{Feller},
\[
\sum_{i > N^{1/u}} \frac{\xi_i^2 i^{-t}}{(1+N i^{-u})^v}
\asymp\sum_{i > N^{1/u}} \S(i)^2 i^{-t-2q-1} \asymp N^{-(t+2q)/u}\S
^2(N^{1/u}).
\]
This is never bigger than the contribution of the first part of the sum,
and of equal order if $(t+2q)/u< v$.
If $(t+2q)/u>v$, then the leading polynomial term
is strictly smaller than $N^{-v}$. If $(t+2q)/u=v$,
then the leading term is equal to $N^{-v}$, but the slowly\vspace*{1pt}
varying part satisfies $\S^2(N^{1/u}) \ll\sum_{i \leq N^{1/u}}{\S
^2(i)}/{i}$,
by Theorem 1(b) on page 281 in \cite{Feller}. Therefore, in both
cases the preceding display is negligible relative to the
first part of the sum.
This proves the final assertion of
the lemma for $c=1$. The proof for general $c>0$ is analogous.
\end{pf}

By the Cauchy--Schwarz inequality, for any $\m\in S^{t/2}$,
\[
\biggl|\sum_i\frac{\xi_i\mu_i}{1+Ni^{-u}}\biggr|^2
\le\|\m\|_{t/2}^2\sum_i \frac{\xi_i^2i^{-t}}{(1+N i^{-u})^2}.
\]
The preceding lemma gives the exact order of the right-hand side.
The application of the Cauchy--Schwarz inequality is sharp, in that
there is equality for some \mbox{$\mu\in S^{t/2}$}. However, this $\mu$ depends
on $N$. For fixed $\mu\in S^{t/2}$ the left-hand side is strictly smaller
than the right-hand side.
%
%le8.3 #&#
%
\begin{lemma}
\label{LemmaTechnicalBias}
For any $t,u\ge0$, $\mu\in S^{t/2}$ and $(\xi_i)$ such that
$|\xi_i|=\break i^{-q-1/2}\S(i)$ for $0<t+2q <2u$ and a slowly varying function
$\S\dvtx (0,\infty)\to(0,\infty)$, as $N \ra\infty$,
\[
\sum_i \frac{|\xi_i\mu_i|}{1+N i^{-u}}
\ll N^{-(t+2q)/(2u)}\S(N^{1/u}).
\]
\end{lemma}
\begin{pf}
We split the series in two parts, and bound the denominator
$1+Ni^{-u}$ by $Ni^{-u}$ or $1$. By the Cauchy--Schwarz inequality,
for any $r>0$,
\begin{eqnarray*}
\biggl|\sum_{i\le N^{1/u}} \frac{|\xi_i\mu_i|}{N i^{-u}}\biggr|^2
&\le&\frac1{N^2}\sum_{i\le N^{1/u}} \frac{\S^2(i)i^r}{i}
\sum_{i\le N^{1/u}}\m_i^2i^{2u-2q-r}\\
&\asymp&\frac1{N^2}\S^2(N^{1/u})N^{r/u}\\
&&{}\times\sum_{i\le N^{1/u}}{\m_i^2i^t}\biggl(\frac{i}{N^{1/u}}
\biggr)^{2u-2q-r-t} N^{(2u-2q-r-t)/u},\\
\biggl|\sum_{i> N^{1/u}} \frac{|\xi_i\mu_i|}{1}\biggr|^2
&\le&\sum_{i> N^{1/u}} \frac{\S^2(i)}{i}i^{-2q}\sum_{i> N^{1/u}}\m_i^2
\asymp\S^2(N^{1/u})N^{-2q/u}\sum_{i> N^{1/u}}\m_i^2.
\end{eqnarray*}
The terms in the remaining series in the right-hand side of the first
inequality are bounded by $\m_i^2i^t$ and tend to zero pointwise
as $N\ra\infty$ if $2u-2q-r-t>0$.
If $t+2q<2u$, then there exists $r>0$ such that the latter
is true, and for this $r$ the sum tends to zero by the dominated
convergence theorem. The other terms collect to
$N^{-(t+2q)/(u)}\S^2(N^{1/u})$. The sum in the right-hand side of the
second inequality is bounded by
$\sum_{i>N^{1/u}}\m_i^2i^t N^{-t/u}=o(N^{-t/u})$.
\end{pf}

%suskaldyti doi

% imsref loaded by lrinkeviciute, 2011-11-11 09:04:07
% imsref loaded by lrinkeviciute, 2011-11-11 09:57:29
%

\printaddresses

\end{document}